\documentclass[10pt,oneside,letterpaper,draft]{amsart}

\usepackage{verbatim}
\usepackage{ifthen}
\usepackage{enumitem}

\usepackage{amsthm}

\usepackage{multirow}

\usepackage[T1]{fontenc}
\usepackage{amsmath}
\usepackage{amssymb}

\usepackage{import}

\usepackage[left=2.75cm,right=2.75cm, top=2.75cm, bottom=2.75cm]{geometry}

\newcommand{\thesisor}[2]{#1}

\usepackage{amsthm}
\usepackage{comment}
\usepackage{enumitem}

\newtheorem{thm}{Theorem}

\newtheorem*{thm*}{Theorem}
\newtheorem{lemma}[thm]{Lemma}
\newtheorem{prop}[thm]{Proposition}
\newtheorem{corr}[thm]{Corollary}
\newtheorem*{corr*}{Corollary}

\theoremstyle{definition}

\numberwithin{thm}{section}
\numberwithin{defn}{section}

\newcommand{\jap}[1]{\langle #1 \rangle }
\newcommand{\floor}[1]{\left\lfloor #1 \right\rfloor }
\newcommand{\R}{\mathbb{R}}
\newcommand{\C}{\mathbb{C}}
\newcommand{\of}{\overline{f}}
\newcommand{\og}{\overline{g}}
\newcommand{\E}{\mathcal{E}}
\newcommand{\F}{\mathcal{F}}
\newcommand{\T}{\mathcal{T}}
\newcommand{\HO}{\mathcal{H}}
\newcommand{\Ham}{\mathcal{H}}
\newcommand{\ft}{\widehat{f}}

\newcommand{\Hs}{{\mathcal{H}^s}}

\renewcommand{\Re}{\text{Re}}
\newcommand{\Sch}{Schr\"{o}dinger}
\renewcommand{\Im}{\mathop{\text{Im}}}
\renewcommand{\Re}{\mathop{\text{Re}}}

\newcommand{\Fmu}{{F_{\mu,\epsilon}}}
\newcommand{\Gmu}{{G_{\mu,\epsilon}}}
\newcommand{\lil}[1]{\frac{\pi}{2}\floor{ \frac{2#1}{\pi} } }
\newcommand{\manys}[2]{\underbrace{#1,\ldots,#1}_{#2\text{ times}}}
\newcommand{\manyf}[1]{\underbrace{f,\ldots,f}_{#1\text{ times}}}

\newcommand{\proj}[1][n]{\Pi_{#1}}

\newcommand{\CS}{Cauchy--Schwarz}

\newcommand{\nrm}[1]{ | #1 |}
\newcommand{\intpi}{\int_{-\pi/4}^{\pi/4}}

\newcommand{\Eqin}{\E_6}
\newcommand{\Tqin}{\T_6}
\newcommand{\Hqin}{\Ham_6}

\newcommand{\Ecub}{\E_4}
\newcommand{\Tcub}{\T_4}
\newcommand{\Hcub}{\Ham_4}

\newcommand{\POR}{}

\renewcommand{\Re}{\text{Re}\,}
\renewcommand{\Im}{\text{Im}\,}

\title[Resonant Hamiltonian systems]{Resonant Hamiltonian systems associated to the one-dimensional nonlinear Schr\"{o}dinger equation with harmonic trapping}
\author{James Fennell}

\address{Courant Institute of Mathematical Sciences, 251 Mercer Street, New York NY 10012, United States}
\email{fennell@cims.nyu.edu}

\numberwithin{equation}{section}
\numberwithin{thm}{section}
\numberwithin{defn}{section}
\numberwithin{example}{section}

\subjclass[2000]{35Q55; 37K05}

\begin{document}

\begin{abstract}

	We study two resonant Hamiltonian systems on the phase space 
        $L^2(\R \rightarrow \C)$:
		the quintic one-dimensional continuous resonant equation,
		and a cubic resonant system that has appeared in 
                the literature as a modified scattering limit
                for an NLS equation with cigar shaped trap.
	We prove that these systems approximate the dynamics of the 
		quintic and cubic one-dimensional NLS with harmonic trapping
		in the small data regime on long times scales.
		We then pursue a thorough study of the dynamics of the resonant systems themselves.
	Our central finding is that these resonant equations 
		fit into a larger class of Hamiltonian systems that have many striking dynamical features:
		non-trivial symmetries such as invariance under the Fourier transform and the flow of the linear \Sch{} equation with harmonic trapping,
		a robust wellposedness theory, including global wellposedness in $L^2$
                and all higher $L^2$ Sobolev spaces, 
		and an infinite family of orthogonal, explicit stationary wave solutions in the form of the Hermite functions.

\end{abstract}

\maketitle 
\renewcommand{\thesisor}[2]{#2}

\section{Introduction}

\noindent
In recent years resonant systems have emerged as extremely useful tools for studying nonlinear \Sch{} equations (NLS).
Resonant equations have been used to construct solutions of the cubic NLS on $\mathbb{T}^2$ that exhibit large growth of Sobolev norms \cite{Colliander2010}.
They have appeared as modified scattering limits for a number of equations, including the cubic NLS on $\R \times \mathbb{T}^d$ \cite{Hani2015B},
	the cubic NLS on $\R^d$ with $2\leq d \leq 5$ and harmonic trapping in all but one direction \cite{Hani2015},
	and a coupled cubic NLS system on $\R \times \mathbb{T}$ \cite{Victor2016}.
The {continuous resonant equation} (CR) was originally 
	shown to approximate the dynamics of small solutions of the two-dimensional cubic NLS on a large torus $\mathbb{T}_L^2$ over
	long times scales (longer than $L^2/\epsilon^2$, where $\epsilon$ is the size of the initial data) \cite{Faou2013}.
	Recent work has extended this by showing that a whole family of CR equations
	approximate the dynamics of NLS on $\mathbb{T}_L^d$ for arbitrary dimension and arbitrary analytic nonlinearity \cite{Buckmaster2017A}.
The original two-dimensional cubic CR equation is the
	same resonant system that appears in the modified scattering limit in \cite{Hani2015} for $d=3$;
	it has also been shown to be a small data approximation for the cubic NLS with harmonic trapping set on $\R^2$ \cite{Germain2016}.

One of the principal reasons that resonant systems are useful is that they generally exhibit a large amount of structure.
They are often Hamiltonian and usually possess many symmetries, a good wellposedness theory, and an infinite number of orthogonal, explicit solutions.
Extensive work has been done on studying such purely dynamical properties of the CR equations:
	starting in the paper that introduced the original two-dimensional cubic equation \cite{Faou2013},
	in subsequent works again on this cubic case \cite{Germain2015,Germain2016},
	and a more recent paper on the general case \cite{Buckmaster2017B}.
This research fits into a larger program of studying the dynamics of nonlocal Hamiltonian PDEs;
	we mention, for example, work on the Szeg\"{o} equation \cite{Gerard2017B} and the lowest Landau level equation \cite{Gerard2017}.

The two-dimensional cubic CR equation has, in particular, been found to have many remarkable dynamical properties.
The PDE is symmetric under many non-trivial actions such as the Fourier transform and the linear flow of the \Sch{} equation (with or without harmonic trapping);
	it is Hamiltonian, and through these symmetries admits a number of conserved quantities.
The equation is globally wellposed in $L^2$ and all higher Sobolev spaces.
It has many explicit stationary wave solutions, including all of the Hermite functions and the function $1/|x|$.
All stationary waves that are in $L^2$ are automatically analytic and exponentially decaying in physical space and Fourier space. 

The present work was initiated by the question of whether these striking properties 
	also hold for the only other continuous resonant equation that scales like $L^2$:
	the one-dimensional quintic continuous resonant equation.
Our investigation subsequently broadened to include another one-dimensional resonant equation
	that is somewhat more physically relevent, and turns out to be the modified scattering limit in \cite{Hani2015} for $d=2$.
Our overall finding is that these Hamiltonian systems do display much of the remarkable dynamical structure of the two-dimensional cubic CR.
In fact, we are able to show that both systems belong to a large
	class of Hamiltonian systems on the phase space $L^2(\R \rightarrow \C)$,
	and that each system is this class bears many of the features of $L^2$ critical CR: 
	    each has a strong symmetry structure, global wellposedness in $L^2$
                and all higher $L^2$ Sobolev spaces,
		and many explicit stationary wave solutions in the form of the Hermite functions.
Typical members of the class lack much of the structure of both cubic two-dimensional and quintic one-dimensional CR -- for example, it is not the case that all $L^2$ stationary waves are analytic -- 
	but 
	our findings do suggest that a number of the properties of the $L^2$ critical CR equations are generic.

\subsection{Presentation of the equations}

The two systems we study in this \thesisor{chapter}{article}
    are resonant systems corresponding
	to
	the nonlinear \Sch{} equation with harmonic trapping,
\begin{equation}
	iu_t - \Delta u + x^2 u = iu_t + Hu = |u|^{2k}u,
	\label{eqn:1:pde}
\end{equation}
where the spatial variable is $x \in \R$ and $k=1,2$ is an integer, so that the nonlinearity is analytic.
The cubic $k=1$ equation is physically relevant: in this case, \eqref{eqn:1:pde} is the {Gross--Pitaevskii} equation and is a model in the physical theory of Bose--Einstein condensates \cite{Gross1961}.

Let us first see how the resonant equations arise.
Looking at the profile $v(t) = e^{-itH}u(t)$ (where $e^{itH}$ is the propagator of the linear equation $iu_t+Hu = 0$),
	we find it satisfies,
        \[
        iv_t = e^{-itH}\left( |e^{itH}v|^{2k} e^{itH} v \right).
        \]
Expressing $v(t)$ in the basis of eigenfunctions of the operator $H$ (namely the Hermite functions), the equation on $v$ can be written as,
\begin{equation}
	iv_t(t) = \sum_{n_1,\ldots,n_{2k+2} \in \mathbb{Z}^+} e^{2itL} \proj[n_{2k+2}]
		\left[ \prod_{m=1}^k \left(  ( \proj[n_m] v(t) )( \overline{ \proj[n_{k+1+m}] v(t) }) \right) \proj[n_{k+1}] v(t) \right],
	\label{eqn:1:pde3}
\end{equation}
where $\proj v$ is the projection onto the eigenspace of $H$ corresponding to eigenvalue $2n+1$.
The phase $L$ in \eqref{eqn:1:pde3} is given by $L = n_1 + \ldots + n_{k+1} -( n_{k+2} + \ldots + n_{2k+2})$.
The resonant terms in the sum in \eqref{eqn:1:pde3} are the terms that are not oscillating in time; that is, those satisfying $L=0$.
The resonant system corresponding to \eqref{eqn:1:pde3} is obtained by considering only the resonant terms.
We will show in Section 2 that this resonant PDE may be written more compactly in terms of a certain time average of the nonlinearity,
\begin{equation}
	iw_t(t) = \frac{2}{\pi} \intpi e^{-isH} \left( | e^{isH} w(t) |^{2k} e^{isH} w(t) \right) ds.
	\label{eqn:1:respde2}
\end{equation}
From this expression we are able to infer that the resonant system
is, up to a rescaling of time,
the Hamiltonian flow on the phase space $L^2(\R \rightarrow \C)$ corresponding to the 
	Hamiltonian,
\begin{equation}
	\Ham_{2k+2}(f) = \frac{2}{\pi} \intpi \int_\R |(e^{itH}f)(x)|^{2k+2} dx dt.
	\label{eqn:1:resham}
\end{equation}

The overall resonant program is to gain information on the dynamics of solutions to \eqref{eqn:1:pde} by studying the associated resonant system \eqref{eqn:1:respde2}.
This program has two, distinct components.
The first is to establish approximation results that rigorously demonstrate that solutions of the resonant system well approximate 
	solutions of the full system in certain function spaces 
        and over certain timescales.
        In Section 2 we prove such an approximation result that is valid for all positive integers $k$.
The second component of the resonant
program is to understand the dynamics of the resonant equation itself.
One then projects these dynamics back to the original equation through the approximation results.
Our analysis focuses on the resonant system \eqref{eqn:1:respde2} 
in depth for the cubic case, when $k=1$, and the quintic case, when $k=2$.
These two cases are particularly significant for separate reasons.
\begin{itemize}
	\item
		The cubic case $k=1$ is physically relevant, as previously mentioned.
		In addition, the resonant equation here is exactly the resonant equation obtained  in \cite{Hani2015} as the modified scattering limit of the NLS equation,
		\begin{equation}
			i u_t - u_{xx} - u_{yy} + |y|^2u = |u|^2u,
			\label{eqn:1:hani}
		\end{equation}
		where the space variable is $(x,y)\in\R^2$.
		Precisely, consider small initial data $u_0(x,y)$.
		Suppose that  $u(x,y,t)$ solves  \eqref{eqn:1:hani} with initial data $(x,y) \mapsto u_0(x,y)$.
		For each fixed $x$, let $w(x,y,t)$ be the solution of the resonant equation \eqref{eqn:1:respde2} with initial data $y \mapsto u_0(x,y)$.
		Then,
		\[
			\lim_{t \rightarrow \infty} \left\| u(x,y,t) - e^{i t (-\partial_{xx} - \partial_{yy} + |y|^2) } w(x,y,2 \ln(t) ) \right\|_{H^N(\R^2)} = 0,
		\]
		where $H^N$ is the usual Sobolev space.
		(This holds for any $N \geq 8$ so long as the initial data is sufficiently small.)
	\item
		In the quintic case, $k=2$, 
		we will prove that the resonant system \eqref{eqn:1:respde2} is precisely the one-dimensional quintic continuous resonant equation.
		It is the only CR equation, other than the original two-dimensional cubic CR equation, that scales like $L^2$.
\end{itemize}

\subsection{Obtained results}

\subsubsection{An approximation theorem}

We begin, in Section 2, by proving the following theorem, which shows that solutions of the resonant equation  \eqref{eqn:1:respde2}
	well-approximate solutions of the full equation \eqref{eqn:1:pde} on a long time scale.
	This theorem is essentially a lower dimensional version of Theorem 3.1
        in \cite{Germain2016}, and our proof follows theirs closely.

	\begin{thm*}[Theorem \ref{thm:2:approx}, page \pageref{thm:2:approx}]
		Define the space $\Hs$ by the norm $\| f \|_{\Hs} = \| H^{s/2} f \|_{L^2}$; this is equivalent to the norm $\| \jap{x}^s f \|_{L^2} + \| \jap{\xi}^s \ft \|_{L^2}$.
	Fix $s>1/2$ and initial data $u_0 \in \Hs$.
	Let $u$ be a solution of the nonlinear \Sch{} equation with harmonic trapping \eqref{eqn:1:pde} 
		and $w$ a solution of the resonant equation \eqref{eqn:1:respde2}, both
		corresponding to the initial data $u_0$.
	Suppose that the bounds $\|u(t)\|_{\Hs}, \| w(t) \|_{\Hs}\leq \epsilon$ hold for all $t\in[0,T]$.
	Then for all $t\in[0,T]$,
	\[
		\| u(t) - e^{itH} w(t)  \|_{\Hs} \leq 
			\left( t (2k+1) \epsilon^{4k+1} + \epsilon^{2k+1} \right) \exp\left( (2k+1) t \epsilon^{2k} \right).
	\]
	In particular if $t \lesssim \epsilon^{-2k}$ then
	$
		\| u(t) - e^{itH} w(t)  \|_{\Hs} \lesssim
			 \epsilon^{2k+1} .
	$
	\end{thm*}

	\subsubsection{Representation formulas for the Hamiltonians}

	Following the approximation result, we focus on 
	studying the resonant system \eqref{eqn:1:respde2} in the cases $k=2$ and $k=1$.
In both cases the rights hand sides of the resonant PDE \eqref{eqn:1:respde2} 
can be written in terms of the multilinear operators,
\[
    \T_{2k+2}(f,\ldots,f) =\frac{2}{\pi} \intpi e^{-isH} ( |e^{isH} f |^{2k}  e^{isH} f ) ds
\]
so that the resonant equation reads $iw_t = \T_{2k+2}(w,\ldots,w)$.
The fact that Hamilton's equation 
can be expressed in terms of multilinear operators
is a nontrivial structural property that
	guides much of the analysis.

In the study of resonant equations, it has turned out to be fundamental 
	to determine alternative representations for the Hamiltonian $\Ham$ 
         and the associated multilinear operator $\T$.
	These alternative representations often reveal structure that is concealed by specific representations such as \eqref{eqn:1:resham}.
	In Sections 3.1 and 4.1 we derive numerous representations for 
        $\Hqin$ and $\Hcub$ respectively.
First, for $\Hqin$, we find the two formulas,
\begin{align}
	\Hqin(f)
	&= \frac{2}{\pi} \int_\R \int_\R  
		| e^{it\Delta}f |^6
		dx dt \label{eqn:1:schham} \\
	&= \frac{1}{\pi^2} \int_{\R^6} 
		f(y_1)
		f(y_2)
		f(y_3)
		f(y_4)
		f(y_5)
		f(y_6) 
		\delta_{y_1+y_2+y_3=y_4+y_5+y_6}
		\delta_{y_1^2+y_2^2+y_3^2=y_4^2+y_5^2+y_6^2} dy, \notag
\end{align}
where, in the first equation, $e^{it\Delta}$ denotes the propagator of the linear \Sch{} equation.
These representations both show that the quintic Hamiltonian system is the one-dimensional quintic continuous resonant equation \cite{Buckmaster2017A}.

To describe our next representations, we require some notation.
For an isometry $A: \R^3 \rightarrow \R^3$, let $E_A$ be the multilinear functional,
\begin{equation}
	E_A(f_1,f_2,f_3,f_4,f_5,f_6) = \int_{\R^3}
						f_1( (Ax)_1 )
						f_2( (Ax)_2 )
						f_3( (Ax)_3 )
						\overline{f_4( x_1 )
						f_5( x_2 )
						f_6( x_3 )} dx_1 dx_2 dx_3,
	\label{eqn:1:EA}
\end{equation}	
where $(Ax)_k = \langle Ax, e_k \rangle$.
The functional $E_A$ is a special case of the type of functional that appears 
on the left hand side in Brascamp-Lieb inequalities \cite{Brascamp1976}.
We then have the following representations:
for the quintic equation, we prove that,
\begin{equation}
	\Hqin(f) = \frac{1}{2\sqrt{3} \pi^2} 
            \int_0^{2\pi} E_{R(\theta)}( f,f,f,f,f,f ) d\theta,
	\label{eqn:1:quinticBL}
\end{equation}	
where $R(\theta)$ is the rotation of $\R^3$ by $\theta$ radians about the axis $(1,1,1)$;
	while for the cubic equation, we prove that,
\begin{equation}
	\Hcub(f) =
        \frac{1}{2\sqrt{2} \pi^2} \int_0^{2\pi} E_{S(\theta)}( G,f,f,G,f,f) d\theta,
	\label{eqn:1:cubicBL}
\end{equation}	
where $G(x) = e^{-x^2/2}$, 
	and $S(\theta)$ is the rotation of $\R^3$ by $\theta$ radians about the axis $(0,1,1)$;

The two representations \eqref{eqn:1:quinticBL} and \eqref{eqn:1:cubicBL} are 
extremely useful for studying $\Hqin$ and $\Hcub$. 
They also place the two Hamiltonians in a larger class of Hamiltonians that, we will find, share much of the same structure.
This is not obvious:
	\emph{a priori} we might expect the Hamiltonians $\Hqin$ and $\Hcub$ to be quite unalike.
The differences in \eqref{eqn:1:quinticBL} and \eqref{eqn:1:cubicBL} are also of note.
The functional $E_A$ has many symmetries, and these
are inherited directly by $\Hqin$.
The presence of the Gaussians $G$ in $\Hcub$
    causes the symmetry group of the cubic equation
	to be smaller.
    It also prevents the cubic Hamiltonian from inheriting the
    scaling law present in $E_A$; 
    this has
	consequences for the possible stationary waves we can construct.

\subsubsection{Properties of the Hamiltonian systems}

Using the representations \eqref{eqn:1:quinticBL}
and \eqref{eqn:1:cubicBL} we begin our analyses of the Hamiltonian 
properties of the resonant equations.

\begin{thm*}[Theorem \ref{thm:4:symE}, page \pageref{thm:4:symE},
and Theorem \ref{thm:5:symE}, page \pageref{thm:5:symE}]
	The Hamiltonians $\Hqin$ and $\Hcub$
        are invariant under the following actions (for any $\lambda$):
	\begin{align*}
		\text{(i) }&\text{Fourier Transform: }f_k \mapsto \ft_k.
		&\text{(vii) }&\text{Linear modulation: }f_k \mapsto e^{i\lambda x  } f_k.\\
		\text{(ii) }&\text{Modulation: }f_k \mapsto e^{i\lambda} f_k.	
		&\text{(viii) }&\text{Translation: }f_k \mapsto f_k(\cdot+\lambda).\\
		\text{(vi) }&\text{\Sch{} with harmonic trapping: }f_k \mapsto e^{i\lambda H} f_k. 
	\end{align*}
        The Hamiltonian $\Hqin$ is, in addition, invariant under the following
        actions (for any $\lambda$):
	\begin{align*}
		\text{(iv) }&\text{Quadratic modulation: }f_k \mapsto e^{i\lambda x^2} f_k.		
		&\text{(v) }&\text{\Sch{} group: }f_k \mapsto e^{i\lambda \Delta} f_k.	\\
	\text{(iii) }&L^2\text{ scaling: }f_k(x) \mapsto \lambda^{1/2} f_k(\lambda x).
	\end{align*}
\end{thm*}

We prove these symmetries by showing that the symmetries actually
hold for the functional $E_A$.
The symmetries are simply inherited by $\Hqin$, while for $\Hcub$
the Gaussian terms in \eqref{eqn:1:cubicBL} prevent the inheritance of some 
symmetries.

These symmetries of $\Hqin$ and $\Hcub$ lead directly,
by Noether's Theorem,
 to 
 conserved quantities for the associated resonant equations.

\begin{corr*}[Table \ref{tbl:4:conserved}, page \pageref{tbl:4:conserved}, and Table \ref{tbl:5:conserved}, page \pageref{tbl:5:conserved}]
	The following are conserved quantities of the resonant equation \eqref{eqn:1:respde2} in the quintic ($k=2$) and cubic ($k=1)$ cases,
	\begin{align*}
		&&\int_\R |f(x)|^2 dx,
		&&\int_\R x|f(x)|^2 dx,
		&&\int_\R if'(x) \of(x)   dx,
		&&\int_\R |x f(x)|^2 + |f'(x)|^2 dx.
	\end{align*}
	In the quintic case $k=2$, we have the additional conserved quantities,
	\begin{align*}
		&&\int_\R \left[ i xf'(x)+f(x) \right] \of(x)  dx,
		&&\int_\R |xf(x)|^2 dx,
		&&\int_\R |f'(x)|^2 dx.
	\end{align*}
\end{corr*}

We next examine the $L^2$ boundedness of the Hamiltonians.

\begin{thm*}[Theorem \ref{thm:4:L2bound}, page \pageref{thm:4:L2bound}, and
    Theorem \ref{thm:5:L2}, page \pageref{thm:5:L2}]
	There hold the sharp bounds,
        \begin{align*}
            \Hqin(f) &= \frac{2}{\pi} \| e^{itH} f \|_{L^6_t L^6_x}^6
                \leq \frac{1}{\sqrt{3\pi}} \| f \|_{L^2}^6,
            &\Hcub(f) &= \frac{2}{\pi} \| e^{itH} f \|_{L^4_t L^4_x}^4
                \leq \frac{1}{\sqrt{2\pi}} \| f \|_{L^2}^4.
        \end{align*} 
        We have equality in the bound for
        $\Hqin$ if and only if $f=\gamma \exp(-\alpha x^2 + \beta x)$
        for $\alpha,\beta,\gamma\in \C$ and $\Re \alpha > 0 $.
        We have equality in the bound for
        $\Hcub$ if and only if $f=\gamma \exp(-(1/2) x^2 + \beta x)$
        for $\beta,\gamma \in \C$.
\end{thm*}

This bound is actually an example of a geometric Brascamp-Lieb inequality, 
	and the classification of the maximizers is already known \cite{Barthe1998}.
We prove the inequality and classify the maximizers
in the present case in a way that appears to be original.

Using the representation $\Hqin(f) = (2/\pi) \| e^{it\Delta}f \|_{{L^6_tL^6_x}}$, from \eqref{eqn:1:schham}, the $L^2$ bound on $\Hqin$ reads,
\[
	\| e^{it\Delta} f \|_{L^6_t L^6_x}^6 \leq \frac{1}{2\sqrt{3}} \|f \|_{L^2}^6,
\]
which is the homogeneous Strichartz inequality in dimension one.
Our work shows that the constant here is the best possible, and that there is equality if and only if $f$ is a Gaussian.
These facts were previously determined in \cite{Foschi2007}.

\subsubsection{Wellposedness of the resonant equations}
We then turn to the PDE problem associated to $\Hqin$ and $\Hcub$;
i.e., the resonant PDE.

\begin{thm*}[Theorem \ref{thm:4:wellposedness}, page \pageref{thm:4:wellposedness},
    and Theorem \ref{thm:5:wellposedness}, page \pageref{thm:5:wellposedness}]
	The mutilinear operators $\Tqin$ and $\Tcub$
        are is bounded from $X^5$ to $X$ for
	(i) $X=L^2$, (ii) $X=L^{2,\sigma}$ for any $\sigma \geq 0$, and (iii) $X=H^{\sigma}$ for any $\sigma \geq 0$.
        (iv) $L^{\infty,s}$, for any $s \geq 1/2$.
        (v) $L^{p,s}$, for any $p \geq 2$ and $s>1/2-1/p$.
\end{thm*}

These bounds lead  directly to local wellposedness  for the resonant equations
in all of these spaces.
By pairing local wellposedness with the conservation of the $L^2$ norm, we get global wellposedness in $L^2$.\POR{}
A persistence of regularity argument then gives global
wellposedness in every $H^\sigma$ for $\sigma>0$.

\begin{corr*}[Theorems \ref{thm:4:wellposedness}, page \pageref{thm:4:wellposedness}, and \ref{thm:5:wellposedness}, page \pageref{thm:5:wellposedness}]
	Hamilton's equations corresponding to $\Hqin$ and $\Hcub$ are
		locally wellposed in $X$ for (i) $X=L^2$, (ii) $X=L^{2,\sigma}$ for any $\sigma \geq 0$, and (iii) $X=H^\sigma$ for any $\sigma \geq 0$.
	They are globally wellposed in $L^2$ and $H^\sigma$.\POR{}
\end{corr*}

It is expected that item (ii) here can be sharpened to show that
$\Tqin$ is bounded from $(\dot{L}^{\infty,1/2})^5$ to $\dot{L}^{\infty,1/2}$
	(here homogeneous weighted $L^\infty$ spaces).
	This is equivalent to $1/\sqrt{|x|}$ being a stationary wave of the quintic Hamiltonian system, which we conjecture.

\subsubsection{Stationary waves}

As with the cubic continuous resonant equation in dimension 2, 
the resonant equations here admit many explicit stationary wave solutions.
\begin{thm*}[Theorem \ref{thm:4:smooth}, page \pageref{thm:4:smooth},
and Theorem \ref{thm:5:smooth}, page \pageref{thm:5:smooth}]
	For every $n \geq 0$, 
        the Hermite function $\phi_n(x)$ is a stationary wave of the Hamiltonian systems $\Hqin$ and $\Hcub$.
\end{thm*}

By letting the symmetries of each of the equations act on $\phi_n$ we can construct more stationary waves; see \eqref{eqn:4:allstats} and \eqref{eqn:5:allstats}.

\begin{thm*}[Theorem \ref{thm:4:smooth}, page \pageref{thm:4:smooth}, and 
Theorem \ref{thm:5:smooth}, page \pageref{thm:5:smooth}
]
	Suppose that $\phi \in L^2$ is a stationary wave solution of the 
        quintic ($k=2)$ or cubic ($k=1$) resonant equation 
	\eqref{eqn:1:respde2}.
	Then there is $\alpha,\beta>0$ such that $\phi e^{\alpha x^2} \in L^\infty$ and $\hat\phi e^{\beta x^2} \in L^\infty$.
	In particular, $\phi$ can be extended to an analytic function on the complex plane.
\end{thm*}

Our proofs of these results rely on a number of refined Strichartz 
estimates \cite{Bernicot2014,Bourgain1998,Tao2003}.
Using the representations for the Hamiltonians
we are able to provide original proofs of these estimates in an elementary way.
For example, in Proposition \ref{thm:4:refined} we prove that
if $\ft_1$ is supported in $B(0,R)^C$ and $\ft_2$ is supported
in $B(0,r)$, for $R>4r$, then,
\[
        \int_{\R^2}
            |e^{it\Delta}f_1|^2
            |e^{it\Delta}f_2|^2
            |e^{it\Delta}f_3|^2
            dx dt
        \leq  \frac{1}{R^{1/6}}
        \| f_1 \|_{L^2}^2
        \| f_2 \|_{L^2}^2
        \| f_3 \|_{L^2}^2.
\]

\subsubsection{Smoothing for the quintic resonant equation}

For the quintic equation, we prove the following smoothing result,
which shows that the operator $\Tqin$ increases the regularity of
Sobolev data.

\begin{thm*}[Theorem \ref{thm:4:smoothing}, page \pageref{thm:4:smoothing}]
For any $\sigma>0$, there is a $\delta>0$ such that
$\Tqin$ is bounded from 
$(L^{2,\sigma})^5$ to $L^{2,\sigma+\delta}$ and
$(H^{\sigma})^5$ to $H^{\sigma+\delta}$.
\end{thm*}

\subsection{Plan of the \thesisor{chapter}{article}}

In Section 2 we prove the main approximation result.
In Section 3 we present results concerning the quintic Hamiltonian system
    defined by $\Hqin$.
Many of our proofs rely on studying the functional $E_A$, and the needed 
properties of $E_A$ are proved as required.
The cubic Hamiltonian system is treated in a similar fashion in Section 4.

\subsection{Notations and conventions}

\begin{itemize}
\item
	For $x\in \R$, the Japanese bracket is $\jap{x} = \sqrt{1+x^2}$.

\item $\jap{f,g}_{L^2} = \int_\R f(x) \og(x) dx$.
	
\item	The Sobolev space $H^\sigma$ is defined by the norm $\|f\|_{H^{\sigma}} = \| \jap{x}^\sigma \ft \|_{L^2}$.
\item	The weighted space $L^{2,\sigma}$ is defined by the norm $\|f\|_{L^{2,\sigma}} = \| \jap{x}^\sigma f \|_{L^2}$.
\item $H=-\Delta + x^2$ is the operator corresponding to the quantum harmonic oscillator.

	\item

		The Fourier transform of $f$ is
$
		\F(f)(\xi)=\hat{f}(\xi) = (2\pi)^{-1/2}  \int_{\R} e^{-ix\xi} f(x) dx.
$
With this convention, the map $f \mapsto \hat{f}$ is an isometry of $L^2(\R)$,
	and the identity $\F(\F(f))(x)=f(-x)$ holds.
We will frequently use the Fourier inversion formula,
		\begin{equation}
			\frac{1}{(2\pi)^n} \int_{\R^n} \int_{\R^n} e^{ia \jap{ w,x}} \phi(w) dw dx
	= \frac{1}{(2\pi)^{n/2}} \int_{\R^n} \hat\phi(aw) dw 
	= \frac{1}{|a|} \phi(0).
			\label{eqn:1:fourierinv}
		\end{equation}

	\item
		We set $G(x)=e^{-x^2/2}$.
For all $a>0$,
$
		\F \left( e^{- \frac{a x^2}{2} }  \right)(\xi) = a^{-1/2} e^{- \frac{ \xi^2 }{ 2a } }
$
and $\int_\R e^{-ax^2} dx = \sqrt{\pi/a}$.

\item $A \lesssim B$ means there is an absolute constant $C$ such that $A \leq CB$. $A \sim B$ means $A \lesssim B$ and $B\lesssim A$.

\end{itemize}

\section{An approximation theorem}
We begin the article by treating more precisely the derivation of the resonant equation \eqref{eqn:1:respde2}
	and then proving the approximation theorem described in the introduction.

Before studying  the nonlinear problem, we recall 
some basic properties of the linear problem corresponding to \eqref{eqn:1:pde}.
These facts will be used extensively throughout the article.
The linear equation corresponding to \eqref{eqn:1:pde} is simply the equation for the quantum harmonic oscillator,
\begin{equation}
	\begin{split}
		i u_t + Hu =
		iu_t - \Delta u + x^2 u &=
		0,	
	\end{split}
	\label{eqn:2:linear}
\end{equation}
where $H = -\Delta + x^2 $.
For any initial data $u_0\in L^2$ there is a unique solution to \eqref{eqn:2:linear}, which we denote $e^{itH} u_0$.
An explicit representation of this solution is given by the \emph{Mehler formula},
\begin{equation}
	e^{itH}u_0(x) = \frac{1}{ \sqrt{ 2\pi |\sin(2t)| } } \int_\R e^{-i \left[ ( x^2/2 + y^2/2)\cos(2t) - xy \right]/\sin(2t) } u_0(y) dy.
	\label{eqn:2:mehler}
\end{equation}
(This and other properties of the linear flow may be found in \cite{Burq2013}.)
From this expression we see that the solution is time-periodic with period $\pi$.

An alternative representation of the solution of \eqref{eqn:2:linear} may be found by examining the \emph{Hermite functions} $\{ \phi_n \}_{n=0}^\infty$.
The Hermite functions are eigenfunctions of $H$ -- they satisfy $H\phi_n = (2n+1)\phi_n$ -- 
	and they form an orthonormal basis of $L^2$.
Each of these functions is a polynomial multiplied by the Gaussian $e^{-x^2/2}$; for example,
	$\phi_0(x) = c_0 e^{-x^2/2}$,
	$\phi_1(x) = c_1 x e^{-x^2/2}$, and
        $\phi_2(x) = c_2 (1-2x^2)e^{-x^2/2}$,
where the constants $c_n$ are normalizing constants that ensure $\| \phi_n \|_{L^2}=1$.
Using the eigenfunction property one finds that $e^{itH} \phi_n = e^{it(2n+1)}\phi_n$.
Let $\proj{} u_0 = \jap{u_0,\phi_n}\phi_n$ be the orthogonal projection onto the eigenspace spanned by $\phi_n$.
Given any $u_0 \in L^2$ we may expand
$
	u_0(x) = \sum_{n=0}^\infty (\proj u_0)(x),
$
and then find,
$
	e^{itH} u_0(x) = \sum_{n=0}^\infty e^{it(2n+1)} (\proj u_0)(x),
$
so the flow has a simple description in the Hermite function coordinates.
We finally note that the Hermite functions satisfy $\phi_n(-x) = (-1)^n \phi_n(x)$;
this may be infered from the formula $\phi_n(x) = c_n e^{x^2/2} (d^n/dx^n) e^{-x^2}$ from \cite{Burq2013}.

We now turn to the nonlinear problem \eqref{eqn:1:pde}. 
The linear part of the equation may be absorbed into the nonlinearity by changing variables to the profile $v(x,t) = e^{-itH}u(x,t)$.
The function $v$ satisfies the equation,
\begin{equation}
	i v_t = e^{-itH} \left( |e^{itH}v|^{2k} e^{itH}  v \right) := N_t(v,\ldots,v),
	\label{eqn:2:vPDE}
\end{equation}
where $N_t$ is the $(2k+1)$ multilinear operator,
\begin{equation}
	N_t( f_1, \ldots, f_{2k+1} ) = e^{-itH}\left[ \left( \prod_{m=1}^k (e^{itH}f_m)( \overline{e^{itH} f_{k+1+m}}) \right) (e^{itH} f_{k+1}) \right].
	\label{eqn:2:Nt}
\end{equation}
We expand each of the functions $f_m$ in the basis of Hermite functions,
$
	e^{itH}f_m = e^{itH} \left(
		\sum_{n_m=0}^\infty \proj[n_m] f_m \right) = 	
		\sum_{n_m=0}^\infty e^{it(2n_m+1)} \proj[n_m] f_m,
$
and then substitute into  \eqref{eqn:2:Nt}.
This yields,
\begin{equation}
	N_t( f_1, \ldots, f_{2k+1} )
		= \sum_{n_1,\ldots,n_{2k+2} \geq 0} 
                \hspace{-0.5cm} e^{2i L t}  \proj[n_{2k+2}] 
		\left[ \left( \prod_{m=1}^k (\proj[n_m] f_m)( \overline{ \proj[n_{k+1+m}] f_{k+1+m}}) \right) \proj[n_{k+1}] f_{k+1} \right],
	\label{eqn:2:Nt2}
\end{equation}
where $L= \sum_{m=1}^{k+1} n_m - n_{k+m+1}$.

In \eqref{eqn:2:Nt2}, when $L\neq0$ the associated term in the sum is oscillating,
	while when $L=0$ the associated term is not.
The resonant equation arises simply from neglecting the oscillatory terms.
Define the multilinear functional $\T$ by
\begin{equation}
	\T(f_1,\ldots,f_{2k+1}) 	
		= \sum_{\substack{n_1,\ldots,n_{2k+2} \geq 0 \\ L=0}}
			\proj[n_{2k+2}] \left[ \left( \prod_{m=1}^k (\proj[n_m] f_m)( \overline{ \proj[n_{k+1+m}] f_{k+1+m}}) \right) \proj[n_{k+1}]  f_{k+1} \right].
	\label{eqn:2:T}
\end{equation}
The resonant PDE is then given by,
\begin{equation}
	iw_t =\T(w,\ldots,w).
	\label{eqn:2:respde}
\end{equation}

\begin{lemma}
	The resonant functional $\T$ is the time average of the functionals $N_r$ over the interval $[-\pi/4,\pi/4]$; that is,
	\begin{equation}
		\T(f_1,\ldots,f_{2k+1}) 	
			= \frac{2}{\pi} \intpi N_r(f_1,\ldots,f_{2k+1}) dr.
		\label{eqn:2:T2}
	\end{equation}
\end{lemma}

\begin{proof}
	We integrate the sum in \eqref{eqn:2:Nt2} over $[-\pi/4,\pi/4]$ term by term.
	If $L=0$ nothing changes and we get the associated term in \eqref{eqn:2:T}.
	If $L$ is even then $\intpi  e^{2iLr} dr = 0$, and the term in \eqref{eqn:2:Nt2} is 0.
	Finally if $L$ is odd, then either $n_{2k+2}$ is even and $L-n_{2k+2}$ is odd, or $n_{2k+2}$ is odd and $L-n_{2k+2}$ is even.
	In the first case we have, using the Hermite function property $(\proj f)(-x) = (-1)^n (\proj f)(x)$, that,
	\begin{align*}
		&\left( \prod_{m=1}^k ((\proj[n_m] f_m)(-x))( \overline{ (\proj[{n_{k+1+m}}] f_{k+1+m})(-x)}) \right) (\proj[{n_{k+1}}] f_{k+1})(-x) 	\\
		& \hspace{2cm}= (-1)^{L-n_{2k+2}}\left( \prod_{m=1}^k (( \proj[{n_m}] f_m)(x))( \overline{ ( \proj[{n_{k+1+m}}] f_{k+1+m})(x)}) \right)
			(\proj[{n_{k+1}}] f_{k+1})(x) ,
	\end{align*}
	and hence the function here is odd. 
	Projecting onto the eigenspace spanned by the even function $\phi_{n_{2k+2}}$ gives the 0 vector.
	The associated term in the sum \eqref{eqn:2:Nt2} is thus 0.
	In the case when $n_{2k+2}$ is even and $L-n_{2k+2}$ is odd a similar analysis shows that the term in the sum is again 0.
	In conclusion, all of terms corresponding to $L\neq 0$ vanish, while those corresponding to $L=0$ are unchanged.
\end{proof}

By virtue of the lemma the resonant equation can be written as,
\begin{align}
	i w_t = \T(w,\ldots,w) 
		= \frac{2}{\pi} \intpi N_r(w(t),\ldots,w(t)) dr    
		&= \frac{2}{\pi} \intpi e^{-irH} \left( |e^{irH} w(t) |^{2k}
                    e^{irH}w(t) \right) dr,
	\label{eqn:2:respde2}
\end{align}
which is precisely \eqref{eqn:1:respde2}.
One can show that the resonant equation is the flow corresponding the Hamiltonian,
$
	\Ham_{2k+2}(f) = \frac{2}{\pi} \intpi \int_\R  | e^{irH} f (x) |^{2k+2} dx dr ,
$
up to a rescaling of time.
The details of this Hamiltonian correspondence are
presented in Theorem \ref{thm:4:hameqn} below.

We now prove the approximation theorem.
The theorem is essentially a lower dimensional analog of Theorem 3.1 in \cite{Germain2016}, and our proof follows theirs closely.
The function space in our theorem is,
\[
	\Hs = \{ u \in L^2 : H^{s/2} u \in L^2 \},
\]
with the norm  $\| u \|_{\Hs} = \| H^{s/2} u \|_{L^2}$.
From \cite{Yajima2008}, we have the norm equivalence 
$
    \| u \|_{\Hs} \sim \| \jap{x}^{s/2}  u \|_{L^2} 
        + \| \jap{\xi}^{s/2} \hat{u} \|_{L^2}.
    $
This space $\Hs$ is useful for two reasons: 
	first, if $s>1/2$, then the space is an algebra (as a direct consequence of the norm equivalence);
	and, second, the space interacts well with the linear propagator $e^{itH}$, as seen in the following Lemma.

\begin{lemma}
	Fix $s \geq 0$.
	For all $u\in \Hs$ and $t \in \R$ we have $\| e^{itH} u \|_{\Hs} \leq \| u \|_{\Hs}$.
	\label{lemma:Hsalgebra}
\end{lemma}

	A general $L^p$ version of this lemma appears in \cite{Bongioanni2006}; for $L^2$, there is the following short proof.

\begin{proof}
	First let $s$ be an even non-negative integer.
        Then $s/2$ is an integer, and it is clear that $H^{s/2}$ 
        commutes with $e^{itH}$.
        By conservation of the $L^2$ norm by $e^{itH}$, we have,
        \[
		\| e^{itH} u \|_{\Hs}
                = \| H^{s/2} e^{itH} u \|_{L^2}
                = \| e^{itH} H^{s/2} u \|_{L^2}
                = \|  H^{s/2} u \|_{L^2}
		= \| e^{itH} u \|_{\Hs}.
        \]
	The result for general $s$ follows from interpolation.
\end{proof}


\begin{thm}
	Fix $s>1/2$ and initial data $u_0 \in \Hs$.
	Let $u$ be a solution of the nonlinear \Sch{} equation with harmonic trapping \eqref{eqn:1:pde} 
		and $w$ a solution of the resonant equation \eqref{eqn:2:respde}, both
		corresponding to the same initial data $u_0$.
	Suppose that the bounds $\|u(t)\|_{\Hs}, \| w(t) \|_{\Hs}\leq \epsilon$ hold for all $t\in[0,T]$.
	Then for all $t\in[0,T]$,
	\[
		\| u(t) - e^{itH} w(t)  \|_{\Hs} \leq 
			\left( t (2k+1) \epsilon^{4k+1} + \epsilon^{2k+1} \right) \exp\left( (2k+1) t \epsilon^{2k} \right).
	\]
	In particular if $t \lesssim \epsilon^{-2k}$ then
	$
		\| u(t) - e^{-itH} w(t)  \|_{\Hs} \lesssim \epsilon^{2k+1} .
	$
	\label{thm:2:approx}
\end{thm}


\begin{proof}
	Let $v(x,t) = e^{-itH}u(x,t)$, so that $v$ satisfies the PDE \eqref{eqn:2:vPDE}.
	We note that $v(x,0) = u(x,0) = u_0(x)$.
	Using the previous Lemma, we find that,
	\begin{equation}
		\| u(t) - e^{itH} w(t)  \|_{\Hs} =
		\| e^{itH} v(t) - e^{itH} w(t)  \|_{\Hs} \leq 
		\| v(t) -  w(t)  \|_{\Hs}.
		\label{eqn:2:uw2vw}
	\end{equation}
	To prove the theorem it therefore suffices to show that $v$ and $w$ are close in $\Hs$.

	Therefore let $v$ and $w$ be solutions of the equations \eqref{eqn:2:vPDE} and \eqref{eqn:2:respde} respectively with the same initial data $u_0$,
	\begin{align}
		i v_t(t) &= N_t(v(t),\ldots,v(t))
			= e^{-itH} \left( |e^{itH} v(t) |^{2k} e^{itH}v(t) \right) , 
                        \label{eqn:2:nonres}	\\
		i w_t(t) &= \T(w(t),\ldots,w(t))	
			= \frac{2}{\pi} \intpi e^{-irH} \left( |e^{irH} w(t) |^{2k} e^{irH}w(t) \right) dr, \label{eqn:2:resthm},
	\end{align}
    and $u_0(x)	= v(x,0)=u(x,0)$.
	Set,
	\begin{align}
		D_t (f_1,\ldots,f_{2k+1}) &= N_t(f_1,\ldots,f_{2k+1})  - \T(f_1,\ldots,f_{2k+1}) 
                \label{eqn:2:Dt}
                \\
		&= \sum_{\substack{n_1,\ldots,n_{2k+2} \geq 0 \\ L \neq 0}} 
                e^{2itL} \proj[{n_{2k+2}}] 
			\left[ \left( \prod_{m=1}^k (\proj[{n_m}] f_m)( \overline{\proj[{n_{k+1+m}}] f_{k+1+m}}) \right) \proj[{n_{k+1}}] f_{k+1} \right].
			\label{eqn:2:Dt2}	
	\end{align}
	From the expressions of the multilinear operators $N_t$ and $\T$ in \eqref{eqn:2:nonres} and \eqref{eqn:2:resthm}
		(or their multilinear versions \eqref{eqn:2:Nt} and \eqref{eqn:2:T2}),
		from Lemma \ref{lemma:Hsalgebra},
		and from the fact that $\Hs$ is an algebra, it follows that $N_t$ and $\T$ are uniformly bounded  from $(\Hs)^{2k+1}$ to $\Hs$.
	The same holds for $D_t$ from \eqref{eqn:2:Dt}.

	Set $\phi(t) = v(t)-w(t)$. 
	Because $\phi(0)=0$, the Duhamel form of the equation on $\phi$ is,
	\begin{equation}
		i \phi(t) = \int_0^t \left[ \T(v(r),\ldots,v(r)) - \T(w(r),\ldots,w(r)) + D_r(v(r),\ldots,v(r)) \right]dr
		\label{eqn:2:duhamel}
	\end{equation}
	We will determine \emph{a priori} bounds on $\phi$.
	For the first term in the integrand here, we can expand by multilinearity to find,
	\begin{align}
		\| \T(v(r),\ldots,v(r)) - \T(w(r),\ldots,w(r)) \|_{\Hs} 
		&\leq \sum_{m=0}^{2k} \| \T(\underbrace{v(r),\ldots,v(r)}_{m \text{ times}},v(r)-w(r),\underbrace{w(r),\ldots,w(r)}_{2k-m\text{ times}} \|_\Hs \notag\\
		& \leq (2k+1) \epsilon^{2k} \| v(r) - w(r)  \|_{\Hs}. \label{eqn:2:thmest1}
	\end{align}

	For the second term in the integrand in \eqref{eqn:2:duhamel} we need to look more closely at the operator $D_t$.
	We first observe the identity,
	$
		e^{2irL} = 
			d/dr \int_{ \lil{r} }^{r} e^{2i \theta L} d\theta,
	$
	where $\floor{x}$ is the smallest integer less that $x$. 
	(Recall from the proof of the first lemma that only even values of $L$ contribute to the sum in \eqref{eqn:2:Dt2}.)
	The interval of integration here has length less than 1.
	We can then handle the second term in \eqref{eqn:2:duhamel} as follows,
	\begin{align*}
		& \int_0^t \left[ D_r(v(r),\ldots,v(r)) \right]ds \\
		&= \sum_{\substack{n_1,\ldots,n_{2k+2} \geq 0 \\ L \neq 0}}  \int_0^t e^{2irL} \proj[ {n_{2k+2}} ]
			\left[ \left( \prod_{m=1}^k (\proj[{n_m}] v(r))( \overline{ \proj[{n_{k+1+m}}] v(r)}) \right) \proj[{n_{k+1}}] v(r)\right] dr\\
		&= \sum_{\substack{n_1,\ldots,n_{2k+2} \geq 0 \\ L \neq 0}}  \int_0^t  
			\frac{d}{dr} \left( \int_{ \lil{r} }^{r} e^{2i\theta L} d\theta \right)
		\proj[{n_{2k+2}}] \left[ \left( \prod_{m=1}^k ( \proj[{n_m}] v)( \overline{ \proj[{n_{k+1+m}}] v}) \right) \proj[{n_{k+1}}] v\right] dr	.
	\end{align*}
	Using integration by parts, we have,
	\begin{align*}
		&\text{(left hand side)}	\\
		&= -\sum_{\substack{n_1,\ldots,n_{2k+2} \geq 0 \\ L \neq 0}}  \int_0^t  
			\left( \int_{ \lil{r}  }^{r} e^{2i \theta L} d\theta \right)
		\frac{d}{dr} \proj[{n_{2k+2}}] \left[ \left( \prod_{m=1}^k (\proj[{n_m}] v  )( \overline{ \proj[{n_{k+1+m}}] v }) \right) \proj[{n_{k+1}}] v\right] dr\\
		& \hspace{2cm}  + \sum_{\substack{n_1,\ldots,n_{2k+2} \geq 0 \\ L \neq 0}}  
			\left( \int_{ \lil{t} }^{t} e^{2i \theta L} d\theta \right)
		\proj[{n_{2k+2}}] \left[ \left( \prod_{m=1}^k (\proj[{n_m}] v)( \overline{ \proj[{n_{k+1+m}}] v}) \right) \proj[{n_{k+1}}] v\right] 	\\
		&= - \int_0^t   \sum_{m=0}^{2k} \int_{\lil{r}}^r D_\theta ( \manys{v(r)}{m},v_r(r), \manys{v(r)}{2k-m} )d\theta dr \\
		& \hspace{2cm}  + 
			\left( \int_{ \lil{t} }^{t} D_\theta(v(t),\ldots,v(t)) d\theta \right).
	\end{align*}
	Because the interval of integration $\left[\lil{t},t\right]$ has length less than 1, we get,
	\begin{equation}
		\begin{split}
			\left\| \int_0^t \left[ D_r(v(r),\ldots,v(r)) \right]ds \right\|_{\Hs} 
				&\leq t (2k+1) \sup_{r \in [0,t]} \left( \| v(r) \|_{\Hs}^{2k} \| v_r(r) \|_{\Hs} \right) + \| v(t) \|_{\Hs}^{2k+1} \\
				&\leq t (2k+1) \epsilon^{4k+1} + \epsilon^{2k+1},
		\end{split}
		\label{eqn:2:thmest2}
	\end{equation}
	where in the last line we have used $\| v_t \|_{\Hs} \leq \| v \|_{\Hs}^{2k+1} \leq \epsilon^{2k+1}$,
		coming from \eqref{eqn:2:vPDE}.

	Combining the estimates \eqref{eqn:2:thmest1} and \eqref{eqn:2:thmest2} we get,
	\[
		\| \phi(t) \|_{\Hs} \leq (2k+1)\epsilon^{2k} 
                \left( \int_0^t \| \phi(s) \|_{\Hs} ds \right) + 
			t (2k+1) \epsilon^{4k+1} + \epsilon^{2k+1}.
	\]
	Gronwell's inequality then implies that,
	\[
		\| v(t) - w(t) \|_{\Hs} = 
		\| \phi(t) \|_{\Hs} \leq 
			\left( t (2k+1) \epsilon^{4k+1} + \epsilon^{2k+1} \right) \exp\left( (2k+1) t \epsilon^{2k} \right),
	\]
	which, with \eqref{eqn:2:uw2vw}, gives the result.
\end{proof}

\section{The quintic resonant equation}
We now turn to the resonant Hamiltonian corresponding the quintic equation,
\begin{equation}
	\Hqin(f) = \frac{2}{\pi} \| e^{itH} f \|_{L^6_t L^6_x }
	= \frac{2}{\pi} \intpi \int_\R | e^{itH}f(x) |^6 dx dt,
	\label{eqn:4:Hhermite}
\end{equation}
which has a corresponding  multilinear functional,
\begin{equation}
	\Eqin(f_1,f_2,f_3,f_4,f_5,f_6) =  \frac{2}{\pi} \intpi \int_\R  
		(e^{itH}f_1)
		(e^{itH}f_2)
		(e^{itH}f_3)
		\overline{(e^{itH}f_4)
		(e^{itH}f_5)
		(e^{itH}f_6)}
		dx dt,
	\label{eqn:4:Ehermite}
\end{equation}
related by $\Hqin(f)=\Eqin(f,f,f,f,f,f)$.
    The functional $\Eqin$ has a large number of permutation symmetries.
	For any two permutations of three elements $\sigma, \sigma' \in S_3$,
	we have,
\begin{equation}
	\Eqin(f_1,f_2,f_3,f_4,f_5,f_6) = 
	\Eqin(f_{\sigma(1)},f_{\sigma(2)},f_{\sigma(3)},f_{\sigma'(4)},f_{\sigma'(5)},f_{\sigma'(6)}) 
	\label{eqn:4:permsym1}
\end{equation}
as well as the symmetry,
\begin{equation}
	\Eqin(f_1,f_2,f_3,f_4,f_5,f_6) = \overline{ \Eqin(f_4,f_5,f_6,f_1,f_2,f_3) }.
	\label{eqn:4:permsym2}
\end{equation}
These symmetries are used to calculate
Hamilton's equation corresponding to $\Hqin$.

\begin{thm}
            \label{thm:4:hameqn}
    Hamilton's equation corresponding to $\Hqin$ is,
    \begin{equation}
        i u_t(t) 
        =\Tqin(u(t),\ldots,u(t)) = \frac{12}{\pi} 
            \intpi  e^{-isH}\left( |e^{isH} u(t) |^4 e^{isH} u(t) \right) ds,
            \label{eqn:4:hameqn}
    \end{equation}
    which is precisely
    the resonant equation \eqref{eqn:1:respde2} up to rescaling of time.
\end{thm}

\begin{proof}
	In order to find Hamilton's equation of motion corresponding to
        $\Hqin$, we first recall the Hamiltonian phase space structure of 
        $L^2(\R \rightarrow \C)$.
	A symplectic form on $L^2$ is given by $\omega(f,g) = -\Im \jap{f,g}_{L^2}$.
	Given a Hamiltonian $\Ham:L^2 \rightarrow \R$, the
        symplectic gradient $\nabla_\omega \Ham$ is defined
        as the unique solution of the equation
		$\omega(\nabla_\omega \Ham(f),g) = 
			\left. d/d\epsilon \right|_{\epsilon=0} \Ham(f+\epsilon g).
                        $
	Hamilton's equation is then $u_t = \nabla_\omega \Ham(u)$.

	In the present case  $\Hqin(f) = \Eqin(f,\ldots,f)$, 
        and we have, by multilinearity,
	\begin{align}
		\left. \frac{d}{d\epsilon} \right|_{\epsilon=0} \Hqin(f+\epsilon g)
		= \sum_{k=1}^{6}  \Eqin( \manyf{k-1},g, \manyf{6-k})
		= 6 \, \Re \, \Eqin( f,f,f,f,f,g  ) ,
		\label{eqn:3:variationR}
	\end{align}
	where in the last step we used the permutation symmetries
        \eqref{eqn:4:permsym1} and \eqref{eqn:4:permsym2}.
	On the other hand, setting,
	$
             i \nabla_w \Hqin(f) = \Tqin(f, \ldots, f),
             $
	we find,
		$\omega(\nabla_\omega \Hqin(f),g) 
			=  -\Im \langle i \Tqin(f, \ldots, f),g \rangle
			=  \Re 
				\langle \Tqin(f, \ldots, f),g \rangle.
	$
        By the definition of the symplectic gradient,
        the right hand sides of this equation and \eqref{eqn:3:variationR}
		must match for all $f$ and $g$.
	By replacing $g$ by $ig$ and using conjugate linearity, we see that this equality condition holding for all $g$ actually implies that,
	\begin{align}
		\langle \Tqin(f,\ldots,f),g \rangle
		= 6 \, \Eqin(f,f,f,f,f,g)
                &= \frac{12}{\pi} \intpi  \int_\R |e^{isH}f|^4 
                    (e^{isH}f)(\overline{e^{isH}g}) dx ds
	\label{eqn:4:TintermsofE}	
                    \\
                &= 
                 \int_\R \left( \frac{12}{\pi} \intpi 
                e^{-isH} \left( |e^{isH}f|^4 (e^{isH}f) \right) \right)
                    \og dx ,
                    \notag
	\end{align}
        where we have used the fact that $e^{isH}$ is an isometry of 
        $L^2(\R \rightarrow \C)$ for all $s$.
        From this equation we determine the formula for $\Tqin(f,\ldots,f)$
        and hence \eqref{eqn:4:hameqn}.
\end{proof}

Our expression of Hamilton's equation is in terms
of the 5-linear map $\Tqin : (L^2)^5 \rightarrow L^2$ defined by duality
in \eqref{eqn:4:TintermsofE}.
It is central to much of the analysis below that 
the equation can be expressed in terms of such a multilinear operator.

Theorem \ref{thm:4:hameqn}
showed that the Hamiltonian flow corresponding to $\Hqin$ is precisely the resonant equation \ref{eqn:2:respde2} in the quintic case $k=2$.
By the approximation result, Theorem \ref{thm:2:approx}, solutions of 
\eqref{eqn:4:hameqn} with initial data of size $\epsilon$ are close to solutions of
$
	i u_t - \Delta u + |x|^2 = |u|^4u
$
in the space $\Hs$ for $s>1/2$ and times $t \lesssim \epsilon^{-5}$.

\subsection{Representations of the Hamiltonian and the flow operator}

A highly useful 
	approach to the study of Hamiltonians such as $\Hqin$ is to determine
	alternative representation formulas for $\Hqin$, $\Eqin$, and $\Tqin$.
	Functionals such as $\Eqin$ can have a large amount of structure that is concealed by a specific representations such as \eqref{eqn:4:Ehermite}.
	This is will be illustrated clearly below.

	First, we show that $\E_6$ is invariant under the Fourier transform.
\begin{lemma}
		\label{thm:4:fourier1}
	The functional $\E_6$ and operator $\T_6$ are invariant under the Fourier transform,
	\begin{align}
		\Eqin(f_1,f_2,f_3,f_4,f_5,f_6) &= \Eqin( \ft_1, \ft_2, \ft_3, \ft_4, \ft_5, \ft_6),
		\label{eqn:4:fourier1}
			\\
		\widehat{\Tqin}(f_1,f_2,f_3,f_4,f_5) 
                &= \Tqin( \ft_1,\ft_2,\ft_3,\ft_4,\ft_5).
		\label{eqn:4:fourier2}
	\end{align}
\end{lemma}

\begin{proof}
	First let $f_k = \phi_{n_k}$ be  Hermite functions.
	Then $\hat\phi_{n_k} = (i)^{n_k} \phi_{n_k}$, and so,
	\begin{equation}
		\Eqin(
			\hat\phi_{n_1},
			\hat\phi_{n_2},
			\hat\phi_{n_3},
			\hat\phi_{n_4},
			\hat\phi_{n_5},
			\hat\phi_{n_6})
		= (i)^{n_1+n_2+n_3-n_4-n_5-n_6} \Eqin(
			\phi_{n_1},
			\phi_{n_2},
			\phi_{n_3},
			\phi_{n_4},
			\phi_{n_5},\phi_{n_6}).	
		\label{eqn:4:hermitefourier}
	\end{equation}
	On the other hand,
        using that $e^{itH}\phi_n = e^{it(2n+1)}\phi_n$, we have,
	\begin{equation}
		\begin{split}
		 \Eqin(\phi_{n_1},\ldots,\phi_{n_6})	
			&= \frac{2}{\pi} \intpi e^{2it(n_1+n_2+n_3-n_4-n_5-n_6)} dt	\\
			&\hspace{2cm}\times\int_\R \phi_{n_1}(x)\phi_{n_1}(x)\phi_{n_2}(x)\phi_{n_3}(x) \overline{ \phi_{n_4}(x)\phi_{n_5}(x)\phi_{n_6}(x) }dx
		\end{split}
			\label{eqn:4:hermiteexpand}
	\end{equation}
	If $n_1+n_2+n_3-n_4-n_5-n_6$ is a nonzero even integer, then the time integral in \eqref{eqn:4:hermiteexpand} is 0.
	If $n_1+n_2+n_3-n_4-n_5-n_6$ is an odd integer, then by the 
        Hermite function property $\phi_{n_k}(-x)=(-1)^{n_k}\phi_n(x)$,
		the integrand in the space integral in \eqref{eqn:4:hermiteexpand} is odd and hence the integral is 0.
	Therefore, using also \eqref{eqn:4:hermitefourier}, if $n_1+n_2+n_3-n_4-n_5-n_6 \neq 0$, both
		$\E(\hat\phi_{n_1},\ldots,\hat\phi_{n_6})$ and 
		$\E(\phi_{n_1},\ldots,\phi_{n_6})$ are 0 and in particular equal.
	Moreover, if $n_1+n_2+n_3-n_4-n_5-n_6=0$, then by \eqref{eqn:4:hermitefourier} 
		$\E(\hat\phi_{n_1},\ldots,\hat\phi_{n_6}) = \E(\phi_{n_1},\ldots,\phi_{n_6})$.

	Because the Hermite functions are a basis of $L^2$, the formula  \eqref{eqn:4:fourier1} holds for all functions $f_k$.
        The statement for $\Tqin$ follows from this and
        \eqref{eqn:4:TintermsofE}.
\end{proof}

\begin{thm}
	There holds the representations,
	\begin{align}
		\Eqin(f_1,f_2,f_3,f_4,f_5,f_6) 
		&= \frac{2}{\pi} \int_\R \int_\R  
			(e^{it \Delta}f_1) 
			(e^{it \Delta}f_2) 
			(e^{it \Delta}f_3) 
			\overline{ (e^{it \Delta}f_4) 
			(e^{it \Delta}f_5) 
			(e^{it \Delta}f_6) }
			dx dt,
		\label{eqn:4:Esch}	\\
		\Tqin(f_1,f_2,f_3,f_4,f_5)(x) &= 
			\frac{12}{\pi} \intpi
			e^{-it\Delta}\left[ (e^{it\Delta}f_1)(e^{it\Delta}f_2)(e^{it\Delta}f_3)
			\overline{(e^{it\Delta}f_4)(e^{it\Delta}f_5)} \right] (x) dt.
	\label{eqn:4:Tsch}
	\end{align}
	\label{thm:4:Esch}
\end{thm}

\begin{proof}
        The lens transform \cite{Tao2009}
        takes solutions $u$ of the linear \Sch{} equation into
		solutions $v$ of the linear \Sch{} equation with harmonic trapping.
	If we let $u_k(x,t)=(e^{itH}f_k)(x)$ and $v_k(x,t) = (e^{it\Delta}f_k)(x)$, 
        the lens transform reads,
	\[
		u_k(x,t)= \frac{ 1 }{\cos(2t)^{1/2}} v_k\left( \frac{ x }{ \cos(2t) }, \frac{ \tan(2t) }{2} \right) e^{i x^2 \tan(2t)/2}.
	\]
	We substitute these expressions into \eqref{eqn:4:Ehermite} and perform two changes of variable.
	In the time variable, we perform $s = \frac{1}{2}\tan(2t)$.
	This change of variables bijectively maps $(-\pi/4,\pi/4)$ to $(-\infty,\infty)$
	and has determinant $\cos(2t)^{-2}$.
	In the space variable we perform $y=x/\cos(t)$; this has determinant $|\cos(2t)|^{-1}$.
	Then,
	\begin{align*}
		\int_{-\pi/4}^{\pi/4} \int_\R (u_1 u_2 u_3\overline{u_4u_5u_6})(x,t) dxdt 
			&=
		\int_0^{\pi/2} \frac{ 1 }{ | \cos(2t)|^3 } \int_\R (v_1v_2v_3\overline{v_4v_5v_6})\left( \frac{x}{\cos(2t)}, \frac{\tan(2t)}{2} \right)(x,t)dxdt 	\\
			&=
		\int_{-\infty}^{\infty} \int_\R (v_1v_2v_3\overline{v_4v_5v_6})\left(y,s  \right)(x,t)dyds ,
	\end{align*}
	which gives \eqref{eqn:4:Esch}.
	The expression for $\Tqin$ follows from this,
	\eqref{eqn:4:TintermsofE},
        and the fact that $e^{it\Delta}$ is an isometry of $L^2$ for
        all $t$.
\end{proof}

\begin{thm}
			\label{thm:4:Tcr}
	Let $\Omega_1(x)=y_1+y_2+y_3-y_4-y_5-x$
	and $\Omega_2(x)=y_1^2+y_2^2+y_3^2-y_4^2-y_5^2-x^2$.
	Then there holds the representations,
	\begin{align}
		\Eqin(f_1,f_2,f_3,f_4,f_5,f_6)
			&=  \frac{1}{\pi^2} \int_{\R^6} f_1(y_1)f_2(y_2)f_3(y_3) \overline{f_4(y_4) f_5(y_5) f_6(y_6)} \delta_{ \Omega_1(y_6) } \delta_{\Omega_2(y_6)} dy,
			\label{eqn:4:Ecr}
				\\
		\Tqin(f_1,f_2,f_3,f_4,f_5)(x)
			&=  \frac{6}{\pi^2} \int_{\R^5} f_1(y_1)f_2(y_2)f_3(y_3) \overline{f_4(y_4) f_5(y_5)}  \delta_{ \Omega_1(x) } \delta_{ \Omega_2(x) } dy.
			\label{eqn:4:Tcr}
	\end{align}
\end{thm}

\begin{proof}
	We evaluate \eqref{eqn:4:Esch} using the 
	fundamental solution formula for the linear \Sch{} equation,
	$(e^{it\Delta}f_k)(x) = (4\pi i t)^{-1/2} \int_\R e^{i|x-y_k|^2/4t} f_k(y_k)dy_k$.
	This gives,
	\begin{align}
		\Eqin&(f_1,f_2,f_3,f_4,f_5,f_6)\\
			&= \frac{2}{\pi} \int_\R \int_\R
				(e^{it\Delta} ft_1 )
				(e^{it\Delta} ft_2 )
				(e^{it\Delta} ft_3 )
				\overline{(e^{it\Delta} ft_4 )
				(e^{it\Delta} ft_5 )
				(e^{it\Delta} ft_6 )} dx dt	\notag \\
			&= \frac{1}{32\pi^4} \int_\R \frac{1}{t^3} \int_\R \int_{\R^6}
				e^{ - i x \Omega_1(y_6)/2t} e^{+i\Omega_2(y_6)/4t }
				f_1(y_1)f_2(y_2)f_3(y_3) \of_4(y_4) \of_5(y_5) \of_6(y_6)
                                dy dx dt	\notag \\
			&= \frac{1}{4\pi^4} \int_\R \int_\R \int_{\R^6}
				e^{+ i x \Omega_1(y_6)} e^{+i s \Omega_2(y_6) }
					 f_1(y_1)f_2(y_2)f_3(y_3) \of_4(y_4) \of_5(y_5) \of_6(y_6) dy dx ds	\notag \\
			&= \frac{1}{2\pi^3} \int_\R \int_{\R^6}
				e^{+i s \Omega_2(y_6) }
					 f_1(y_1)f_2(y_2)f_3(y_3) \of_4(y_4) \of_5(y_5) \of_6(y_6) \delta_{\Omega_1(y_6)} dy dx	\label{eqn:4:thestart}\\
			&= \frac{1}{\pi^2}  \int_{\R^6}
					 f_1(y_1)f_2(y_2)f_3(y_3) \of_4(y_4) \of_5(y_5) \of_6(y_6) \delta_{\Omega_1(y_6)} \delta_{\Omega_2(y_6)} dy, \notag
	\end{align}
	which is \eqref{eqn:4:Ecr}.
	Equation \eqref{eqn:4:Tcr} follows immediately from definition \eqref{eqn:4:TintermsofE}
		with the $L^2$ inner product integration in $y_6$.
\end{proof}

\begin{thm}
	\label{thm:4:EBL1}
	There holds the representations,
	\begin{align}
	\Eqin(f_1,f_2,f_3,f_4,f_5,f_6)
		&= \frac{1}{2\pi^2}  \int_{\R^6}
		f_1(\beta + \xi) f_2(\lambda \beta + \gamma) f_3( \lambda \gamma + \xi - \lambda \xi) \notag 	\\
		&\hspace{2.2cm}\of_4(\lambda \beta + \xi) \of_5(\beta + \lambda \gamma + \xi - \lambda \xi) \of_6(\gamma) 
	 d\beta d\eta d\xi d\gamma  
		\label{eqn:4:EBL1}.	\\
	\Tqin(f_1,f_2,f_3,f_4,f_5)(x) 
		&= \frac{3}{\pi^2}  \int_{\R^6}
		f_1(\beta + \xi) f_2(\lambda \beta + x) f_3( \lambda x + \xi - \lambda \xi) 	\notag \\
		&\hspace{2.2cm}\of_4(\lambda \beta + \xi) \of_5(\beta + \lambda x+ \xi - \lambda \xi) 
	 d\beta d\eta d\xi 
		\label{eqn:4:TBL1}.
	\end{align}
\end{thm}

\begin{proof}
	We start with formula \eqref{eqn:4:thestart}.
	Introduce new variables $\alpha,\beta,\gamma,\eta,\xi$
	by $y_1 = \beta+\xi$, $y_2=\eta+\gamma$, $y_3=\alpha$, $y_4=\eta+\xi$ and $y_5=\alpha+\beta$.
	We calculate $y_6 = y_1+y_2+y_3-y_4-y_5 = \gamma$ and 
	\[
		\Omega_2(y_6)= y_1^2+y_2^2 + y_3^2 - y_4^2 - y_5^2 - y_6^2 = 2\beta \xi + 2\gamma\eta - 2 \eta \xi - 2\alpha \beta,
	\]
	which gives the formula,
	\begin{align*}
	\Eqin(f_1,f_2,f_3,f_4,f_5,f_6)
	= \frac{1}{2\pi^3}  \int_{\R^6}
		&e^{2it\left[ \beta\xi + \gamma\eta - \eta \xi - \alpha \beta \right]}	
		f_1(\beta+\xi) f_2(\eta+\gamma) f_3(\alpha)  \\
		&\of_4(\eta+\xi) \of_5(\alpha+\beta) \of_6(\gamma) 
	d\alpha d\beta d \gamma  d\eta d\xi dt.
	\end{align*}
	Now change variables from $\eta$ to $\lambda$ through $\eta=\lambda\beta$.
	This gives $d\eta = |\beta| d\lambda$ and therefore,
	\begin{align*}
	\Eqin(f_1,f_2,f_3,f_4,f_5,f_6)
	= \frac{1}{2\pi^3}  \int_{\R^6}
		&|\beta|e^{2it\beta \left[  \xi +  \gamma \lambda  - \xi\lambda - \alpha \right]}	
		f_1(\beta+\xi) f_2(\lambda \beta + \gamma) f_3(\alpha) \\
		&\of_4(\lambda \beta+\xi) \of_5(\alpha+\beta) \of_6(\gamma) 
	d\alpha d\beta d \gamma  d\eta d\xi dt.
	\end{align*}
	Next we use the Fourier inversion formula 
$
	\int_\R\int_\R e^{iatx}\phi(x)dxdt = 
			2\pi |a|^{-1} \phi(0),
$
with $a=2\beta$  and $x(\alpha)=   \xi +  \gamma \lambda  - \xi\lambda - \alpha  $.
This gives,
	\begin{align*}
	\E_A(f_1,f_2,f_3,f_4,f_5,f_6)
	= \frac{1}{2\pi^2}  \int_{\R^6}
		f_1&(\beta+\xi) f_2(\lambda \beta+\gamma) f_3(\xi + \lambda \gamma - \lambda \xi) 	\\
		&\of_4(\lambda \beta+\xi) \of_5(\beta + \xi + \lambda \gamma - \lambda \xi) \of_6(\gamma) 
	 d\beta d\eta d\xi d\gamma  ,
	\end{align*}
	which is \eqref{eqn:4:EBL1}.
		The representation \eqref{eqn:4:TBL1} follows from the definition of $\Tqin$ in \eqref{eqn:4:TintermsofE}
		with the $L^2$ inner product integration in $\gamma$.
\end{proof}

Before starting the next result we recall the following notation from the 
introduction:
for any matrix $A: \R^3 \rightarrow \R^3$ the multilinear functional $E_A$ 
is defined by,
\begin{equation}
    E_A(f_1,f_2,f_3,f_4,f_5,f_6) = \int_{\R^3} 
        f_1((Ax)_1) f_2((Ax)_2) f_3((Ax)_3) \overline{
            f_4( x_1 ) f_5( x_2 ) f_6( x_3 ) 
        }
            dx_1 dx_2 dx_3,
            \label{eqn:4:EAdefn}
\end{equation}
and $T_A$ is defined by duality using the formula,
\begin{equation}
    \langle T_A(f_1,\ldots,f_5),g \rangle = E_A(f_1,\ldots, f_5, g).
    \label{eqn:4:TA}
\end{equation}

\begin{thm}
	\label{thm:4:EBL2}
	There holds the representations,
	\begin{align}
		\Eqin(f_1,f_2,f_3,f_4,f_5,f_6) &=  \frac{1}{2\pi^2} \int_\R \frac{1}{ \lambda^2 - \lambda + 1} E_{A(\lambda)}(f_1,f_2,f_3,f_4,f_5,f_6) d\lambda,
		\label{eqn:4:EBL2}	\\
		\Tqin(f_1,f_2,f_3,f_4,f_5)(x) &= \frac{1}{2\pi^2} \int_\R \frac{1}{ \lambda^2 - \lambda + 1} T_{A(\lambda)}(f_1,f_2,f_3,f_4,f_5)(x) d\lambda,
		\label{eqn:4:TBL2}
	\end{align}
	where, for all $\lambda$, $A(\lambda)$ is an isometry and $A(\lambda)(1,1,1)=(1,1,1)$.
	(The matrix $A(\lambda)$ is given explicitly in \eqref{eqn:4:alambda1} below.)
\end{thm}

\begin{proof}
	In formula \eqref{eqn:4:EBL1}, let $y_1,y_2,y_3$ be the arguments of $f_1,f_2,f_3$ respectively,
		and let $x_1,x_2,x_3$ be the arguments of $\of_4,\of_5,\of_6$ respectively.
	We have,
	\[	
		y = \begin{pmatrix}
			y_1 \\ y_2 \\ y_3
		\end{pmatrix}
		=\begin{pmatrix}
			\beta + \xi 	\\	\lambda\beta + \gamma	\\ \lambda \gamma + \xi - \lambda \xi
		\end{pmatrix}
		=\begin{pmatrix}
			1	&0	&1	\\
			\lambda	&1	&0	\\
			0	&\lambda	&1-\lambda
		\end{pmatrix}
		\begin{pmatrix}
			\beta \\ \gamma \\ \xi
		\end{pmatrix}
		:= B(\lambda)
		\begin{pmatrix}
			\beta \\ \gamma \\ \xi
		\end{pmatrix},
	\]
	and,
	\[	
		x = \begin{pmatrix}
			x_1 \\ x_2 \\ x_3
		\end{pmatrix}
		=\begin{pmatrix}
			\lambda \beta + \xi \\ \beta + \lambda \gamma + \xi - \lambda \xi \\ \gamma
		\end{pmatrix}
		=\begin{pmatrix}
			\lambda	&0	&1	\\
			1	&\lambda	&1-\lambda	\\
			0	&1	&0
		\end{pmatrix}
		\begin{pmatrix}
			\beta \\ \gamma \\ \xi
		\end{pmatrix}
		:= C(\lambda)
		\begin{pmatrix}
			\beta \\ \gamma \\ \xi
		\end{pmatrix}.
	\]
	In equation \eqref{eqn:4:EBL1}, perform the linear change of variables $x = C(\lambda)(\beta,\gamma,\xi)$.
	We find that $\det C(\lambda) = \lambda^2 - \lambda + 1 = (\lambda-\frac{1}{2})^2 + \frac{3}{4}>0$; 
		in particular $C(\lambda)^{-1}$ is defined for all $\lambda$.
	Let $A(\lambda) = B(\lambda)C(\lambda)^{-1}$.
	Changing variables then establishes
		\eqref{eqn:4:EBL2}.
		The expression for $\Tqin$ follows using the definition of $T_A$ 
                \eqref{eqn:4:TA}.

	A calculation reveals that,
	\begin{align}
		A(\lambda) &=
		B(\lambda) C(\lambda)^{-1} 
		=
			\frac{1}{\lambda^2-\lambda+1} 
		\begin{pmatrix}
			\lambda	&1-\lambda	&\lambda^2-\lambda	\\
			\lambda^2-\lambda	&\lambda	&1- \lambda	\\
			1-\lambda	&\lambda^2-\lambda	&\lambda
		\end{pmatrix}.\label{eqn:4:alambda1}
	\end{align}
	It remains to verify the two properties of $A(\lambda)$.
	These can, of course, be determined from the formula \eqref{eqn:4:alambda1};
		however it is more insightful to see how they arise naturally from the 
		combinatorical structure of the arguments to the functions in \eqref{eqn:4:EBL1}.
		\begin{enumerate}[label=(\roman*)]
			\item
				By inspecting \eqref{eqn:4:EBL1}, we find that the squares of the arguments in $f_1,f_2,f_3$ sum to the squares of the arguments in $f_4,f_5,f_6$,
				\begin{equation}
					(\beta+\xi)^2 + (\lambda\beta + \gamma)^2 + (\lambda\gamma+\xi-\lambda\xi)^2
					=(\lambda\beta)^2 +(\beta + \lambda\gamma + \xi - \lambda \xi)^2 + (\gamma)^2.
					\label{eqn:3:sumsquares}
				\end{equation}
				This gives, for all $x \in \R^3$, that $| B(\lambda)x |^2 = \sum_{k=1}^3 | \jap{B(\lambda),e_k}|^2 = \sum_{k=1}^3 | \jap{C(\lambda),e_k}|^2 = |C(\lambda)x|^2$.
				Setting $x = C(\lambda)^{-1}y$ gives $|A(\lambda)y|^2 = |y|^2$ for all $y\in\R^3$, and hence $A(\lambda)$ is an isometry.
			\item
				Again in \eqref{eqn:4:EBL1}, we see that the arguments in $f_1,f_2,f_3$ sum to the arguments in $f_4,f_5,f_6$,
				\begin{equation}
					(\beta+\xi) + (\lambda\beta + \gamma) + (\lambda\gamma+\xi-\lambda\xi)
					=(\lambda\beta)+(\beta + \lambda\gamma + \xi - \lambda \xi) + (\gamma).
					\label{eqn:3:sumargs}
				\end{equation}
				Setting $e=(1,1,1)$, this means that for all $x$, $\jap{B(\lambda)x,e}=\jap{C(\lambda)x,e}$.
				Set $y=C(\lambda)x$ to give $\jap{A(\lambda)y, e}=\jap{y,e}$.
				Because $A$ is an isometry, $A^* = A^{-1}$, and so $\jap{y,A^{-1}e}=\jap{y,e}$ for all $y$, and hence $Ae=e$.
		\end{enumerate}
		We note that the expressions \eqref{eqn:3:sumsquares} and \eqref{eqn:3:sumargs} arise naturally from the $\delta$ arguments in \eqref{eqn:4:Ecr}.
		The properties of $A(\lambda)$ in (i) and (ii) should 
                be considered generic for continuous resonant type equations.
\end{proof}

\begin{thm}
	\label{thm:4:EBL3}
	We have the representations,
	\begin{align}
		\E(f_1,f_2,f_3,f_4,f_5,f_6) &=  \frac{1}{2 \sqrt{3} \pi^2} \int_{0}^{2\pi} E_{R(\theta)}(f_1,f_2,f_3,f_4,f_5,f_6) d\theta,	
		\label{eqn:4:EBL3}		\\
		\T(f_1,f_2,f_3,f_4,f_5)(x) &= \frac{\sqrt{3}}{\pi^2} \int_0^{2\pi} T_{R(\theta)}(f_1,f_2,f_3,f_4,f_5)(x) d\theta,
		\label{eqn:4:TBL3}
	\end{align}
	where $R(\theta)$ is the rotation of $\theta$ radians about the axis $(1,1,1)$.
\end{thm}

\renewcommand{\comment}[1]{}

\begin{proof}
	Because the matrix $A(\lambda)$ is an isometry, $\det A(\lambda)=+1$, and $A(\lambda)(1,1,1) = (1,1,1)$,
	the matrix must, in fact, be a rotation about the axis $(1,1,1)$.
	For any rotation $A$ of $\R^3$, the angle of rotation $\theta$ satisfies, $2\cos(\theta) + 1 = \mathop{\text{Trace}}(A)$.
	In the present case, this means,
	\begin{equation}
		\cos(\theta) =\phi(\lambda):= \frac{1}{2} \left( \text{Trace}(A(\lambda)) - 1 \right) = \frac{1}{2} \left( \frac{3\lambda}{ \lambda^2-\lambda+1}  - 1 \right).
		\label{eqn:4:t2l}
	\end{equation}
	The formula \eqref{eqn:4:EBL3} follows from performing the
        bijective change of variables $\lambda \mapsto\theta$.
\end{proof}

        \comment{
        , which we do rather carefully.

	From analyzing \eqref{eqn:4:t2l}, we determine that $\phi$ has the following properties.
	It satisfies $\phi(-1)=-1$, $\phi(1)=1$; $\phi$ is increasing on $[-1,1]$; $\phi$ is decreasing on $(-\infty,-1] \cup [1,\infty)$;
		and $ \lim_{\lambda \rightarrow -\infty} \phi(\lambda) = \lim_{\lambda \rightarrow +\infty} \phi(\lambda) = 0$.
	By setting $\lambda = 0$ in \eqref{eqn:4:alambda2}, we determine that $\sin(\theta)=-1\sqrt{3}<0$, and hence $\theta = -4\pi/3 \in [\pi,2\pi]$.
	On the other hand, as $\lambda \rightarrow \pm \infty$, $\sin(\theta) \rightarrow + 1/\sqrt{3}>0 $.
	From these considerations and continuity, we infer that 
		that under $\lambda \mapsto \theta$, $[-1,1]$ is bijectively mapped to $[\pi,2\pi]$,
		while $(-\infty,-1)\cup(1,\infty)$ is bijectively mapped to $ (0,\pi)$.
	In all, $\R$ is bijectively mapped to $(0,2\pi]$.

	To perform the change of variables, we need to find the determinent which is given by,
	\[
		\frac{ d\theta }{ d\lambda } =\left| \frac{d}{d\lambda} \arccos \left( \frac{1}{2} \left[ \frac{ 3\lambda }{ \lambda^2 - \lambda + 1 } - 1 \right] \right) \right|.
	\]
	To simplify the computation, we first find that if $a=\sqrt{3}(1-\lambda)/(1+\lambda)$, then
	\begin{align*}
		\arccos \left( \frac{1}{2} \left[ \frac{ 3\lambda }{ \lambda^2 - \lambda + 1 } - 1 \right] \right) 
		&=
		\arccos \left( \frac{ 1 - a^2 }{ 1+ a^2 } \right) = \arctan \left( \frac{ 2a }{ 1-a^2 } \right) 	\\
		&	= 2 \arctan\left( a \right) 
			= 2 \arctan\left( \frac{ \sqrt{3}(1-\lambda)}{ 1 +\lambda  }\right) ;
	\end{align*}
	we then calculate $|d\theta/d\lambda| = \sqrt{3} (\lambda^2-\lambda+1)^{-1}$.
	Formula \eqref{eqn:4:EBL3} then follows.
        }

\subsection{Symmetries of the Hamiltonian and conserved quantities of the flow}

\begin{thm}
	\label{thm:4:symE}
	The functional $\Eqin(f_1,f_2,f_3,f_4,f_5,f_6)$ is invariant under the following actions (for all $\lambda$).
	\begin{enumerate}[label=(\roman*)]
		\item	Fourier transform, $f_k \mapsto \ft_k$.
		\item	Modulation, $f_k \mapsto e^{i\lambda} f_k$.
		\item	$L^2$ scaling, $f_k(x) \mapsto \lambda^{1/2} f_k(\lambda x)$.
		\item	Linear modulation, $f_k \mapsto e^{i\lambda} f_k$.
		\item	Translation, $f_k \mapsto f_k(\cdot+\lambda)$.
		\item	Quadratic modulation, $f_k \mapsto e^{i\lambda x^2} f_k$.
		\item	\Sch{} group, $f_k \mapsto e^{i\lambda \Delta} f_k$.
		\item	\Sch{} with harmonic trapping group, $f_k \mapsto e^{i\lambda H} f_k$.
	\end{enumerate}
\end{thm}

\begin{proof}
We will prove that if a matrix $A: \R^3 \rightarrow \R^3$
is an isometry and satisfies $A(1,1,1) = (1,1,1)$ then 
the functional $E_A$
as defined in \eqref{eqn:4:EAdefn} is
invariant under all of these symmetries.
By the representation \eqref{eqn:4:EBL3} for $\Eqin$,
these symmetries are inherited by $\Eqin$.

(i) Because $A$ is an isometry, we have $\langle \xi, Ax \rangle = \langle A^{-1} \xi, x \rangle$
	for all $\xi,x \in \R^3$.
Now calculating,
\begin{align*}
	E_A(\ft_1,\ldots,\ft_{6})	
	&= \frac{1}{(2\pi)^3} \int_{\R^3}  \prod_{k=1}^3
					\left( \int_{\R} e^{-i\xi_k (Ax)_k } f_k(\xi_k) d\xi_k \right)
					\left( \int_{\R} e^{i\nu_k x_k } \of_{3+k}(\nu_k) d\nu_k \right)
				dx		\\
		&
		= \frac{1}{(2\pi)^3}   \int_{\R^3} \int_{\R^{2\times 3}} e^{ -i \langle A^{-1}\xi -  \nu, x \rangle }
			\prod_{k=1}^3 f_k(\xi_k) \of_{3+k}(\nu_k) d\xi d\nu dx,
\end{align*}
	where in the last line we have used $\jap{\xi,Ax}+\jap{\nu,x}=\jap{A^{-1}\xi-\nu,x}$.
	We first change variables $y(\xi)= A^{-1}\xi - \nu$, or $\xi(y) = Ay+ A\nu$.
The determinant of this change of variables is 1 because $A$ is an isometry.
Performing the change of variables then gives the required identity,
	\begin{align*}
	E_A(\ft_1,\ldots,\ft_{6})	
		&= \frac{1}{(2\pi)^3}   \int_{\R^3} \int_{\R^{2 \times 3}} e^{ -i \langle y, x \rangle }
			\prod_{k=1}^3 f_k(A\nu_k + Ay) \of_{3+k}(\nu_k) dy d\nu dx	\\
			&=
		  \int_{\R^3} 
			\prod_{k=1}^3 f_k((A\nu)_k) \of_{3+k}(\nu_k) d\nu 
                        = E_A(f_1,\ldots,f_{6}),
	\end{align*}
	where in the second equality we used the Fourier inversion identity \eqref{eqn:1:fourierinv} with $a=1$.

(ii) This is clear from the definition of $E_A$.

(iii) Let $f_k^\lambda(x) = \lambda^{1/2} f_k(\lambda x)$.
	Writing out $E_A$ and performing the change of variables 
		$y = \lambda x$ (with $dy = \lambda^3 dx$) 
                gives the relation 
                $
		E_A(f_1^\lambda,\ldots,f_{6}^\lambda) 
			=  E_A(f_1,\ldots,f_{3}).$

    (iv) Because $A$ is an isometry,
		$\nrm{Ax}^2=\nrm{x}^2$ for all $x \in \R^3$.
	Using this, we have,
	\begin{align*}
		E_A( e^{i\lambda |x|^2} f_1,\ldots, e^{i\lambda |x|^2} f_{6} ) 
			&= \int_{ \R^3 } \prod_{k=1}^3 e^{ i\lambda |(Ax)_k|^2 } f_k( (Ax)_k ) e^{-i\lambda |x_k|^2} \of_{3+k}( x_k ) dx 	\\
			&= \int_{ \R^3 } e^{i \lambda |Ax|^2} e^{-i\lambda |x|^2 } \prod_{k=1}^n f_k( (Ax)_k ) \of_{3+k}( x_k ) dx 
			= E_A(f_1,\ldots,f_{6}).
	\end{align*}

    (v) Using the previous part and the invariance of $E_A$
    under the Fourier transform from part (a), we find,
		$E_A( e^{i\lambda \Delta} f_1,\ldots, e^{i\lambda \Delta} f_{6} ) 
                =
		E_A( e^{-i\lambda |x|^2} \ft_1,\ldots, e^{-i\lambda |x|^2} \ft_{6} ) 
		    = E_A(\ft_1,\ldots,\ft_{6})    
		    = E_A(f_1,\ldots,f_{6}).$

(vi)	In this part we use $t$ instead of $\lambda$, and show invariance of the functional under $e^{itH}$.
	First, we note that if $n$ is an integer then $e^{i(\pi/2+n\pi)H} f = \ft$ (from, for instance, the Mehler formula \eqref{eqn:2:mehler}).
	The $t = \pi/2+n\pi$ case thus follows from part (i).
	If $t \neq \pi/2+n\pi$ then we may again represent $e^{itH}f$ 
        in terms of $e^{is\Delta}g$ using the lens transform.
	There holds,
		\begin{equation}
		(e^{itH}f_k)(x) =
			\frac{1}{ \sqrt{ \cos(2t) } } (e^{i(\tan(2t)/2)\Delta}f_k) \left( \frac{ x}{\cos(2t)} \right) e^{ix^2 \tan(2t)/2} .
			\label{eqn:3:lenstransform}
		\end{equation}
	We substitute this expression into the functional.
		Using in turn the symmetries (iv) (with $\lambda = \tan(2t)/2$),
                (iii) (with $\lambda = 1/\cos(2t)$), and 
                (v) (with $\lambda = \tan(2t)/2)$, we determine that,
	\begin{align*}
		E_A( e^{ i t H } f_1, \ldots, e^{i t H} f_{6} ) 
			&=
			E_A \left(
				\frac{1}{ \sqrt{ \cos(2t) } } (e^{i(\tan(2t)/2)\Delta}f_1) \left( \frac{ x}{\cos(2t)} \right),\ldots
                                \right)
                                \\
			&=
			E_A \left(
				(e^{i(\tan(2t)/2)\Delta}f_1) \left( x \right),\ldots,
				(e^{i(\tan(2t)/2)\Delta}f_{6}) \left( x  \right) \right) = E_A( f_1,\ldots,f_{6}).
	\end{align*}

(vii)
Let $e=(1,1,1) \in \R^3$. We have,
	\begin{align*}
		E_A( e^{i\lambda x} f_1,\ldots, e^{i\lambda x}f_{2n} ) 
			&= \int_{ \R^n } \prod_{k=1}^n e^{ i\lambda (Ax)_k } f_k( (Ax)_k ) e^{-i\lambda x_k} \of_{n+k}( x_k ) dx	\\
			&= \int_{ \R^n } e^{i\lambda \langle Ax, e \rangle } e^{-\lambda \langle x, e \rangle} \prod_{k=1}^n  f_k( (Ax)_k ) \of_{n+k}( x_k ) dx
                        = E_A(f_1,\ldots,f_{2n}),
	\end{align*}
	where in the last step we used $\langle Ax,e\rangle= \langle x, A^{-1}e \rangle=\langle x,e\rangle$.

(viii)
	This follows immediately from the previous part and the invariance of the functional  under the Fourier transform, as in item (iv),
		noting that the Fourier transform takes
		$x\mapsto e^{i\lambda x  } f(x)$ to 
		$\xi \mapsto \ft(\xi+\lambda)$.
\end{proof}

\begin{corr}\label{thm:4:TAcomm1} 
	We have the following commuter equalities,
	\begin{align}
		e^{i\lambda Q} \Tqin(f_1,f_2,f_3,f_4,f_5) &= 
                \Tqin(
                    e^{i\lambda Q} f_1,
                    e^{i\lambda Q} f_2,
                    e^{i\lambda Q} f_3,
                    e^{i\lambda Q} f_4,
                    e^{i\lambda Q} f_5
                )	,	\label{eqn:4:TAcomm1} 
	\end{align}
	\begin{align}
		Q \Tqin(f_1,f_2,f_3,f_4,f_5) &= \Tqin(Qf_1,f_2,f_3,f_4,f_5) + 
			\Tqin(f_1,Qf_2,f_3,f_4,f_5) 
                        + \Tqin(f_1,f_2,Qf_3,f_4,f_5) \notag \\ 
			&\hspace{1cm}	- \Tqin(f_1,f_2,f_3,Qf_4,f_5) - \Tqin(f_1,f_2,f_3,f_4,Qf_5),
                \label{eqn:4:TAcomm2} 
	\end{align}
	where $Q$ are the operators: $Q=1$, $Q=x$, $Q=id/dx$, $Q=x^2$, $Q=\Delta$, $Q=H$.
\end{corr}

\begin{proof}
	For each of the operators $Q$, the flow map $e^{i\lambda Q}$ is an isometry of $L^2$ for all $\lambda$, and,
	\[
		\Eqin( e^{i\lambda Q} f_1,\ldots, e^{i\lambda Q} f_{6} ) 
		= \Eqin( f_1,\ldots, f_{6} ) ,
	\]
	from Theorem \ref{thm:4:symE}.
	For each $g \in L^2$, we thus have,
	\begin{align*}
		\langle e^{i\lambda Q} \Tqin(f_1,\ldots,f_{5}),g \rangle_{L^2} 
			&= \langle  \Tqin(f_1,\ldots,f_{5}),e^{-i\lambda Q} g \rangle_{L^2}
        = \Eqin(f_1,\ldots,f_5, e^{-i \lambda Q}g)
        \\&
        = \Eqin( e^{i\lambda Q} f_1,\ldots, e^{i\lambda Q} f_5,g)
	= \langle \Tqin(e^{i \lambda Q} f_1,\ldots,e^{i \lambda Q} f_{2n-1}), g\rangle_{L^2},
	\end{align*}
	which gives \eqref{eqn:4:TAcomm1}.
	To get \eqref{eqn:4:TAcomm2}, differentiate \eqref{eqn:4:TAcomm1} with respect to $\lambda$ and set $\lambda=0$. 
\end{proof}

Because each of the flows $e^{itQ}$ can be realized as a Hamiltonian flow,
Noether's Theorem gives that
the Hamiltonian flow associated to $\Hqin$ has conserved quantities 
associated to symmetries (ii) through (viii).
These symmetries and conserved quantities and summarized in Table \ref{tbl:4:conserved}.

\begin{table}
        \caption{Symmetries of $\Hqin$ and conserved quantities of the 
        quintic resonant equation.}
				\label{tbl:4:conserved}
{\centering

	\renewcommand{\arraystretch}{1.5}
	\begin{tabular}{|ccc|}\hline
		Symmetry of $\Hqin$		&Conserved quantity	&Operator commuting with $\Tqin$\\	\hline
				$f\mapsto e^{i\lambda}f$	&$\int_\R |f(x)|^2 dx$		&$1$	\\
				$f\mapsto f_\lambda$	&$\int_\R \left[ i xf'(x)+f(x) \right] \of(x)  dx$ &	\\
				$f\mapsto e^{i\lambda x}f$	&$\int_\R x|f(x)|^2 dx$	&$x$	\\
				$f\mapsto f(\cdot+\lambda)$	&$\Re\int_\R f'(x) \of(x)   dx$	 &$id/dx$\\
				$f\mapsto e^{i\lambda |x|^2}f$	&$\int_\R |xf(x)|^2 dx$	 &$x^2$	\\
				$f\mapsto e^{i\lambda \Delta}f$	&$\int_\R |f'(x)|^2 dx$	 &$\Delta$\\
				$f\mapsto e^{i\lambda \HO}f$	&$\int_\R |x f(x)|^2 + |f'(x)|^2 dx$ &$H$	\\ \hline
	\end{tabular}

        }

\end{table}

\subsection{Boundedness of the functional and wellposedness of Hamilton's equation}

\begin{thm}
	There holds the following sharp bound,
	\begin{equation}
		|\Eqin(f_1,f_2,f_3,f_4,f_5,f_6)| \leq \frac{1}{\pi\sqrt{3}} 
			\prod_{k=1}^6 \| f_k \|_{L^2};
			\label{eqn:4:bound}
	\end{equation}
	which means in particular $0 \leq \Hqin(f) \leq 1/(\pi\sqrt{3}) \| f \|_{L^2}^6$.
	Equality holds in \eqref{eqn:4:bound} if and only if each $f_k$ is the same Gaussian $\gamma e^{-\alpha x^2 + \beta x}$ for some $\alpha,\beta, \gamma \in \C$ 
        and $\Re \alpha >0$.
	\label{thm:4:L2bound}
\end{thm}

\begin{proof}
First we let $A:\R^3 \rightarrow \R^3$ be a linear isometry.
We have,
\begin{align}
    |E_A(f_1,f_2,f_3,f_4,f_5,f_6)|
    &\leq 
    \int_{\R^3} 
        |f_1(x_1)f_2(x_2)f_3(x_3)| \cdot
        |f_4((Ax)_1)f_5((Ax)_2)f_6((Ax)_3)| dx_1 dx_2 dx_3  \notag \\
    &\leq
    \left( \int_{\R^3} 
        |f_1(x_1)f_2(x_2)f_3(x_3)|^2 dx_1 dx_2 dx_3 \right)^{1/2} \notag \\
    &\hspace{1.5cm}
    \left( \int_{\R^3} 
        |f_4((Ax)_1)f_5((Ax)_2)f_6((Ax)_3)|^2 dx_1 dx_2 dx_3 \right)^{1/2}.
        \label{eqn:4:EAcs}
\end{align}
In the second integral we perform the change of variables $y=Ax$.
The change of variables has determinent 1, because $A$ is an isometry,
and hence the second term is transformed into a term identical in structure
to the first.
In the first term and the second term we can
integrate over
each variable seperately, and hence determine that
$
    |E_A(f_1,f_2,f_3,f_4,f_5,f_6)| \leq \| f_1 \|_{L^2} \cdots \|f_6\|_{L^2}.
$

Now turning to $\Eqin$, using 
	representation \eqref{eqn:4:EBL3} we have,
	\begin{equation}
		|\Eqin(f_1,f_2,f_3,f_4,f_5,f_6)| 
			\leq \frac{1}{2 \sqrt{3} \pi^2} \int_0^{2\pi} | E_{R(\theta)}(f_1,f_2,f_3,f_4,f_5,f_6)|d\theta 
			 \leq \frac{1}{\pi \sqrt{3}} 
			 \prod_{k=1}^6 \| f_k \|_{L^2},
                \label{eqn:4:L2bnd}
	\end{equation}
	which is the inequality \eqref{eqn:4:bound}.

	For equality to hold, we must have equality in \eqref{eqn:4:EAcs}
        for almost every $R(\theta)$; namely we must have 
        $
		|E_{R(\theta)}(f_1,f_2,f_3,f_4,f_5,f_6)| = 
			\prod_{k=1}^6 \| f_k \|_{L^2}
			\label{eqn:4:AlambdaAE}
        $
	for almost every $\theta\in[0,2\pi]$.
        Our use of the Cauchy-Schwartz inequality
        in \eqref{eqn:4:EAcs} means that this happens if and 
        only if,
	\begin{equation}
		f_1( (R(\theta) x)_1 ) 
		f_2( (R(\theta) x)_2 ) 
		f_3( (R(\theta) x)_3 ) =
		f_{4}( x_1 )
		f_{5}( x_2 )
		f_{6}( x_3 ),
		\label{eqn:4:cseq}
	\end{equation}
        for almost every $\theta\in[0,2\pi]$ and $(x_1,x_2,x_3) \in \R^3$.
        Using the fact that $R(\theta)$ is an isometry and 
        that $R(\theta)(1,1,1) = R(\theta)$ one readily verifies
        that this equality does hold
        if the functions $f_k$ are the same Gaussian $e^{- \alpha x^2 +\beta x}$ for $\alpha,\beta \in \C$ and $\Re \alpha>0$.
        The inequality \eqref{eqn:4:L2bnd} is thus sharp.

The converse statement, that functions $f_1,\ldots,f_6$ satisfy
\eqref{eqn:4:cseq} only if each of the functions $f_k$
is the same Gaussian is more involved.
A proof specifically adapted to the present circumstance is presented in
\cite{MyThesis}.
However the $L^2$ equality on $E_A$
is in fact a special case of a geometric Brasscamp-Lieb inequality \cite{Ball2006},
and hence the inequality being saturated by Gaussians is a special
case of the general theorem in \cite{Barthe1998}.
\end{proof}

\begin{thm}
	We have the operator bound
	$
		\| \Tqin(f_1,f_2,f_3,f_4,f_5) \|_X \leq C_X \prod_{k=1}^5 \| f_k \|_X,
	$
	for the following spaces.
	\begin{enumerate}[label=(\roman*)]
		\item $X=L^2$ with $C_X=2\sqrt{3}/\pi$,
		\item $X=L^{2,\sigma}$, for any $\sigma \geq 0$.
		\item $X=H^{\sigma}$, for any $\sigma \geq 0$.
		\item $X=L^{\infty,s}$, for any $s>1/2$.
		\item $X=L^{p,s}$, for any $p\geq 2$ and $s > 1/2 - 1/p$.
	\end{enumerate}
	\label{thm:4:Tbounds}
\end{thm}

\begin{proof}
(i) Follows by duality \eqref{eqn:4:TintermsofE} and
the $L^2$ bound on $\Eqin$ \eqref{eqn:4:bound}.

(ii) Let $\langle x \rangle = \sqrt{1+x^2}$ be the Japanese bracket.
We will first show that
$\Eqin(f_1,\ldots, f_5 , \langle t \rangle^{\sigma} g) 
\leq \Eqin( \langle t \rangle^\sigma f_1, \ldots, \langle t \rangle^\sigma f_5, g)$,
and then determine the bound on $\Tqin$ by duality.

Let $A: \R^3 \rightarrow \R^3$ be an isometry.
Fix $x \in \R^3$.
	Because $A$ is an isometry we have,
	$
		|x_3|^2 \leq \nrm{x}^2 = \nrm{Ax}^2 = \sum_{k=1}^3 |(Ax)_k|^2.
	$
	Therefore there is an integer $l$ such that $|x_3|^2 \leq 3 | (Ax)_l |^2$.
        With $\langle\rangle$ denoting the Japanese bracket,
        we then have $\langle x_3 \rangle \leq 3 \jap{ (Ax)_l }$ and so,
    $
		\jap{ x_3 } \leq 3 \left( \prod_{k=1,2}^n \jap{ (Ax)_k) } \jap{ x_k } \right) \jap{ (Ax)_3},
	$	because in all cases $\jap{t} \geq 1$.
	In terms of the functional $E_A$, this gives,
	\begin{equation}
			E_A(|f_1|,
                            \ldots, |f_5|, \jap{t}^\sigma |f_{6}| ) 
			\leq 3^\sigma	E_A(\jap{t}^\sigma |f_1|,\ldots, 
                        \jap{t}^\sigma |f_{5}|, |f_6| ),
		\label{eqn:3:weighttrans}
	\end{equation}
        and the same inequality is inherited by $\Eqin$ by \eqref{eqn:4:EBL3}.

	Now applying this to $\Tqin$,  we have,
	\begin{align*}
		\langle \Tqin(f_1,\ldots,f_{2n-1}), g \rangle_{L^{2,\sigma}} 
		    &=
		\langle 
                      \Tqin(f_1,\ldots,f_{2n-1}), \jap{t}^{2\sigma}  g \rangle_{L^{2}}	
                      =6 \, \Eqin(f_1,\ldots,f_5,  \jap{t}^{2\sigma}   g)    \\
                & \leq 
                      6\cdot 3^\sigma \, \Eqin( \jap{t}^{\sigma}   f_1,\ldots,
                       \jap{t}^{\sigma}    f_5,  \jap{t}^{\sigma}   g)    
	    < \left( 20 \cdot 3^{\sigma} \prod_{k=1}^{2n-1} \| f_{k} \|_{L^{2,\sigma}} \right)
				\| g \|_{L^{2,\sigma}} ,
	\end{align*}
	which gives the result for $X=L^{2,\sigma}$.

(iii) This follows from (ii) and using the invariance of the operator
$\Tqin$ under the Fourier transform as given in Theorem \ref{thm:4:fourier1}.

	The bound (iv) is proved in Theorem \ref{thm:4:Linftys} below.

	The bound (v) comes from interpolating between the bounds in (iv) and (ii).
\end{proof}

\begin{thm}
	\label{thm:4:wellposedness}
	Consider the Cauchy problem,
	\begin{equation}
		\begin{split}
			iu_t &= \Tqin(u,u,u,u,u),		\\
			u(t=0) &= u_0,
		\end{split}
		\label{eqn:4:cauchy2}
	\end{equation}
	which is Hamilton's equation corresponding to $\Hqin$ and the resonant equation \eqref{eqn:1:respde2} in the quintic $k=2$ case up to rescaling by time.
	
	\begin{enumerate}[label=(\roman*)]
		\item The Cauchy problem \eqref{eqn:4:cauchy2} 
                is locally wellposed in $X$ for any of the
                spaces X in Theorem \ref{thm:4:Tbounds}.
		\item The Cauchy problem \eqref{eqn:4:cauchy2}  
                is globally wellposed in $L^2$.\POR{}
		\item Propagation of regularity:
                the equation is globaly wellposed in $H^\sigma$ for every $\sigma>0$.
	\end{enumerate}
\end{thm}

\begin{proof}
        (i) The Duhamel formulation of the Cauchy problem \eqref{eqn:4:cauchy2} is,
        \[
            u(t) = R\left[ u(t) \right]
            = u_0 + \int_0^t \Tqin( u(s),u(s),u(s),u(s),u(s) ) ds
        \]
        Using multilineariy of $\Tqin$ and the bounds in Theorem \ref{thm:4:Tbounds}
        it is easy to show that for any $T>0$ 
        there is an $\epsilon$ ball around $0$ in
        the space $C_0([0,T],X)$ on which
        $R$ is a contraction mapping.
        Local wellposedness then follows from Banach's Fixed Point Theorem.

                (ii)
                By Banach's Fixed Point Theorem,
                        the local time of existence of a solution to \eqref{eqn:4:cauchy2} in $L^2$ depends only on $\| u_0 \|_{L^2}$.
			Because $\| u \|_{L^2}$ is conserved by the flow \eqref{eqn:4:cauchy2}, 
			by the usual argument the $L^2$ solution is global.

        (iii) This is classical.
        We have $(d/dt) \| u \|_{H^\sigma} \lesssim \| u \|_{H^\sigma} \| u\|_{L^2}^2$.
        From this we see that the $H^\sigma$ norm cannot blow up,
        and hence the $H^\sigma$ solution is global.
\end{proof}

\subsection{Analysis of the stationary waves}

Stationary wave solutions are solutions of the form $e^{i\omega t } \psi(x)$ for some $\omega \in \R$ and a function $\psi$.
By substitution into \eqref{eqn:4:hameqn}, we find that $\psi$ must satisfy,
\begin{equation}
	- \omega \psi(x) = \Tqin(\psi,\psi,\psi,\psi,\psi)(x).
	\label{eqn:4:stationary}
\end{equation}

\begin{thm}
The Hermite functions are stationary waves.
That is, for all $n\geq 0$ there is a number $\omega_n \in \R$
such that $u(x,t) = e^{it \omega_n} \phi_n(x)$
is an explicit solution of \eqref{eqn:4:hameqn}.
\end{thm}

\begin{proof}
    In the proof of Lemma \ref{thm:4:fourier1} we found that for indices
    $n_1,\ldots,n_6$ we had $\Eqin(\phi_{n_1},\ldots, \phi_{n_6})=0$
    unless $n_1+n_2+n_3=n_4+n_5+n_6$.
    In particular, setting $n=n_1=\cdots=n_5$ and $m=n_6$, we have
    \[  
        \langle
            \Tqin(\phi_n,\ldots,\phi_n), \phi_m
        \rangle =
            \Eqin(\phi_n,\ldots,\phi_n, \phi_m) = 0,
    \]  
    unless $n=m$.
    Because the Hermite functions are a basis of $L^2$
    and 
        $ \Tqin(\phi_n,\ldots,\phi_n) \in L^2$ by Theorem \ref{thm:4:Tbounds},
    this implies that , 
            $\Tqin(\phi_n,\ldots,\phi_n) = \omega_n \phi_n$
            for some $\omega_n \in \R$.
    The result follows.
\end{proof}

By letting the symmetries of $\Tqin$ act on $\phi_n$, we find that each of the functions
\begin{equation}
	a e^{i bx + icx^2} \phi_n( dx + e ),
	\label{eqn:4:allstats}
\end{equation}
for $a \in \C$ and $b,c,d,e\in \R$ is a stationary wave solution of
\eqref{eqn:4:hameqn}.

\subsubsection{Regularity of stationary waves: introduction}

All of the stationary waves \eqref{eqn:4:allstats} are analytic and decay in space like $e^{-\alpha x^2}$ for some $\alpha \in \R$.
The remainder of this section is devoted to a proof \emph{any} function $\psi\in L^2$ satisfying \eqref{eqn:4:stationary} 
	is automatically analytic and exponentially decaying in space like $e^{-\alpha x^2}$.
Our proof follows closely the proof of the analogous result for the two-dimensional continuous resonant equation in \cite{Germain2016},
	which in turn is based on work in \cite{Hundertmark2009};
	there are also similar results in \cite{Erdogan2011,Green2016}.

Our proof here has two main ingredients.
Roughly speaking, once a multilinear functional can supply these ingredients,
the associated Hamiltonian system will
	satisfy a result like Theorem \ref{thm:4:strongL2} below. 
The first ingredient is an ability to transfer exponential weight from
	one input of the functional to the other inputs.
The second ingredient is a refined multilinear estimate.

\subsubsection{Exponential weight transfer}

For fixed $\mu,\epsilon>0$, define,
\[
	G_{\mu, \epsilon}(x) = \exp\left( \frac{ \mu x^2 }{ 1 + \epsilon x^2} \right).
\]

\begin{lemma}
	\label{thm:3:weight}
	If $\{ f_1,\ldots,f_6 \}$ are positive functions, then
	\begin{equation}
		\Eqin( f_1,\ldots, f_{5}, f_{6} \Gmu )
		\leq 
                \Eqin( f_1\Gmu ,\ldots, f_{5}\Gmu , f_{6} ).
                \label{eqn:4:weight} 
	\end{equation}
\end{lemma}

\begin{proof}
We will prove the result for the functional $E_A$ where $A$ is an isometry;
the result for $\Eqin$ then follows from the representation \eqref{eqn:4:EBL3}.

	Define $F_{\mu,\epsilon} = \mu |x|/(1+\epsilon |x|)$,
        so that $ \Gmu(x) = \exp(\Fmu(x^2))$.
	We record two properties of $\Fmu$.
	First, for $x>0$, $\Fmu$ is increasing,
        as may be seen from a simple calculation of the derivative.
	Next, we have $\Fmu(x_1+x_2) \leq \Fmu(x_1)+\Fmu(x_2)$.
	This may be seen from,
	\begin{align*}
		\Fmu(x_1+x_2) &= \Fmu( |x_1+x_2|) 	
			\leq \Fmu( |x_1| + |x_2|)		
			= \mu \frac{ |x_1| + |x_2| }{ 1 + \epsilon |x_1| + \epsilon |x_2| }	\\
			&\leq \mu \frac{ |x_1| }{ 1 + \epsilon |x_1| }
						+ \mu \frac{ |x_2| }{ 1 + \epsilon |x_2| } = \Fmu(x_1) + \Fmu(x_2).
	\end{align*}

	Now because $A$ is an isometry we have, for all $x\in \R^3$,
	$
		x_3^2 = \sum_{k=1}^3 (Ax)_k^2 - \sum_{k=1}^{2} (x_k)^2
	$
	and hence by the sublinearity property of $\Fmu$,
	\begin{align*}
		\Fmu( x_3^2	)	=
		 \Fmu \left( \sum_{k=1}^3 (Ax)_k^2 + \sum_{k=1}^{2}- (x_k)^2 \right)	
		\leq \sum_{k=1}^3 \Fmu( (Ax)_k^2 ) + \sum_{k=1}^{2} \Fmu(x_k^2).
	\end{align*}
	Then, because $x \mapsto e^x$ is increasing,
	$
		\Gmu(x_3) = \exp( \Fmu(x_3^2) ) \leq
			\prod_{k=1}^3 \Gmu( (Ax)_k )
			\prod_{k=1}^{2} \Gmu( x_k  ).
	$
	Applying this to $E_A$, we have,
	\begin{align}
		E_A(f_1,\ldots,f_{5},f_{6} \Gmu ) &= \int_{\R^n} \left(
		 \prod_{k=1}^{2} f_k( (Ax)_k )\of_{3+k} ( x_k ) \right) 
                 f_3( (Ax)_3) \of_{6}(x_3) {\Gmu(x_3) } dx \notag	\\
			&\leq \int_{\R^3} \left(
		 \prod_{k=1}^{2} f_k( (Ax)_k ) \Gmu( (Ax)_k )	
                 \of_{3+k} ( x_k ) \Gmu(x_k)  \right)
                \notag \\&\hspace{4cm} f_3( (Ax)_3) {\Gmu((Ax)_3)} \of_{6}(x_3) dx 	
                 \notag \\
		&=
		E_A(f_1 \Gmu,\ldots,f_{5} \Gmu,f_{6}  ),
                \label{eqn:4:EAweight}
	\end{align}
	which is what we wanted to prove.
\end{proof}

\subsubsection{Refined multilinear Strichartz estimates}

The second ingredient we need is a so-called refined multilinear Strichartz estimate.
Such estimates are treated in a number of works \cite{Bernicot2014,Bourgain1998,Tao2003}.
Lemma 111 in \cite{Bourgain1998} is prototypical of the type of estimate we require here:
	it states that if functions $f_1,f_2\in L^2( \R^2 \rightarrow \C)$ satisfy supp$\ft_1 \subset B(0,N)$ and 
	supp$\ft_2 \subset B(0,M)^C$, with $N \ll M$, then,
	\[
		\| (e^{it\Delta}f_1) (e^{it\Delta} f_2) \|_{L^4(\R^2 \times \R)} \lesssim \left( \frac{N}{M} \right)^{1/2} 
		\| f_1 \|_{L^2(\R^2)} \| f_2 \|_{L^2(\R^2)} .
	\]
	The right hand side is decaying for large $M$ and small $N$.
In our case, 
	under similar support assumptions on functions $\ft_i$ and $\ft_j$,
	we would like to have analogous  control on,
	\[
		\Eqin(f_1,f_2,f_3,f_4,f_5,f_6) 
		= \frac{2}{\pi} \int_\R \int_\R  
			(e^{it \Delta}f_1) 
			(e^{it \Delta}f_2) 
			(e^{it \Delta}f_3) 
			\overline{ (e^{it \Delta}f_4) 
			(e^{it \Delta}f_5) 
			(e^{it \Delta}f_6) }
			dx dt;
	\]
	namely, we would like an $L^2$ bound that is decaying as the supports of $\ft_i$ and $\ft_j$ become further and further apart.
	Using the representations \eqref{eqn:4:EBL1} and \eqref{eqn:4:EBL2} we are in fact
	able to determine the required refined multilinear estimate in an elementary way.

	Because we know that $\Eqin$ is invariant under the Fourier transform, it is equivalent to state
		the support assumptions in terms of $f_i$ and $f_j$ and not their Fourier transforms.

\begin{prop}
	Suppose that the support of $f_2$ is in $B(0,R)^C$ 
	and the supports of $f_3$, $f_5$ and $f_6$ are in $B(0,r)$, with $R>4r$.
	Then
	\[
		|\Eqin(f_1,f_2,f_3,f_4,f_5,f_6)|	\leq \frac{1}{R}
			 \prod_{k=1}^6 \| f_k \|_{L^2}.
	\]
\end{prop}

\begin{proof}
We use the representation of $\Eqin$ given in \eqref{eqn:4:EBL1},
	\begin{align*}
	\Eqin(f_1,f_2,f_3,f_4,f_5,f_6)
		&= \frac{1}{2\pi^2}  \int_{\R^6}
		f_1(\beta + \xi) f_2(\lambda \beta + \gamma) f_3( \lambda \gamma + \xi - \lambda \xi) \notag 	\\
		&\hspace{2.2cm}\of_4(\lambda \beta + \xi) \of_5(\beta + \lambda \gamma + \xi - \lambda \xi) \of_6(\gamma) 
	 d\beta d\eta d\xi d\gamma  .
	\end{align*}
	We identify
	a large set in $\lambda$ on which the integrand is 0.
	We will then use the representation \eqref{eqn:4:EBL2} to obtain $L^2$ bounds,
	recalling that the integrand as a function of $\lambda$ is the same in both representations.

Under the assumptions of the proposition, the integrand is non-zero only when
$
	|\beta| \leq |\beta + \lambda \gamma + \xi - \lambda \xi | + |\lambda\gamma + \xi - \lambda \xi| \leq 2r,
$
and only when,
$
	|\lambda \beta| \geq |\lambda \beta + \gamma | - |\gamma| 
        \geq R-2r \geq R/2.
        $
It follows that the integrand is non-zero only when $|\lambda|> R/4$.
	Then, using the representation \eqref{eqn:4:EBL2},
	\begin{align*}
		|\Eqin(f_1,f_2,f_3,f_4,f_5,f_6) | 
		&\leq
		\frac{1}{2\pi^2} \int_{ |\lambda|>R/4 } \frac{1}{\lambda^2-\lambda+1} |E_{A(\lambda)}(f_1,f_2,f_3,f_4,f_5,f_6)| d\lambda	\\
		&\leq \frac{1}{\pi^2} \left( \int_{ |\lambda|>R/4 } \frac{1}{\lambda^2-\lambda+1} d\lambda  \right)
			 \prod_{k=1}^6 \| f_k \|_{L^2}
		\leq  \frac{1}{R} 
			 \prod_{k=1}^6 \| f_k \|_{L^2},
	\end{align*}
        using the $L^2$ bound on $E_A$ from Theorem \ref{thm:4:L2bound}.
\end{proof}

\begin{prop}
			 \label{thm:4:refined}
	Suppose that for some $i$ and some $j$, the support of $f_i$ is in $B(0,R)^C$ 
		and the support of $f_j$ is in $B(0,r)$, with $R>4r$.
	Then
	\begin{equation}
		|\Eqin(f_1,f_2,f_3,f_4,f_5,f_6)|	\leq \frac{1}{R^{1/6}}
			 \prod_{k=1}^6 \| f_k \|_{L^2}.
			 \label{eqn:4:refined}
	\end{equation}
\end{prop}

\begin{proof}
	We assume, by rescaling, that $\| f_k \|_{L^2} = 1$ for all $k$.
	We have the crude bound $\Hqin(f) = \| e^{it\Delta}f_k \|_{L^6}^6 \leq 1$.
	Then,
	\begin{align*}
		| \Eqin(f_1,f_2,f_3,f_4,f_5,f_6)|
			& \leq \frac{2}{\pi} \int_{\R^2} |e^{it\Delta}f_1| \cdot |e^{it\Delta}f_2| \cdot |e^{it\Delta}f_3| \cdot |e^{it\Delta}f_4| \cdot |e^{it\Delta}f_5| \cdot |e^{it\Delta}f_6 |dxdt,	\\
			&= \frac{2}{\pi} \left\| (e^{it\Delta}f_1) (e^{it\Delta}f_2) (e^{it\Delta}f_3) (e^{it\Delta}f_4) (e^{it\Delta}f_5) (e^{it\Delta}f_6) \right\|_{L^1}	\\
			&\leq \frac{2}{\pi} \| (e^{it\Delta}f_i) (e^{it\Delta}f_j) \|_{L^3}	
			= \frac{2}{\pi} \| (e^{it\Delta}f_i)^3 (e^{it\Delta}f_j)^3 \|_{L^1}^{1/3}	\\
			&\leq \frac{2}{\pi} \| (e^{it\Delta}f_i) (e^{it\Delta}f_j)^2\|_{L^2}^{1/3}	 \| (e^{it\Delta}f_i)^2 (e^{it\Delta}f_j) \|_{L^2}^{1/3}	\\
			&\leq \left( \frac{2}{\pi}  \right)^{2/3} \Eqin(f_j,f_i,f_j,f_i,f_j,f_j)^{1/6}  \leq \frac{1}{R^{1/6}},
	\end{align*}
	which is \eqref{eqn:4:refined}.
\end{proof}

\subsubsection{Regularity of stationary waves}

Using the weight transfer property \eqref{eqn:4:weight} and the refined multilinear Strichartz estimate \eqref{eqn:4:refined},
	 we prove that stationary waves are necessarily analytic.
We begin with an integrability result.

\begin{thm}
	\label{thm:4:strongL2}
	Suppose $\phi \in L^2$ satisfies
	\begin{equation}
		|\omega| |\phi(x)| \leq \T (|\phi|,|\phi|,|\phi|,|\phi|,|\phi|)(x).
		\label{eqn:4:weakSW}
	\end{equation}
	Then there exists $\alpha>0$ such that $x \mapsto \phi(x) e^{\alpha x^2} \in L^2$.
\end{thm}

\begin{proof}[Proof of Theorem \ref{thm:4:strongL2}]
	For the proof, we will find $\mu$ so that we have the bound
		$\| \phi G_{\mu,\epsilon} \|_{L^2} \lesssim 1$
	independently of $\epsilon$. 
	Taking the limit $\epsilon \rightarrow 0$ will the yield the result.

	We can clearly assume that $\phi(x) \geq 0$, and will do so throughout.
	For any $M>0$ define,
	\begin{align*}
		\phi_<(x) &= \phi(x) \chi_{|x|\leq M}(x),
		&\phi_{\sim}(x) &= \phi(x) \chi_{M<|x|\leq M^2}(x),
		&\phi_>(x) &= \phi(x) \chi_{M^2<|x|}(x).
	\end{align*}
	We have the decomposition $\phi = \phi_> + \phi_{\sim} + \phi_>$, and the supports are all disjoint, which gives,
	\[
		\| \phi \Gmu \|_{L^2}^2 =
		\| \phi_< \Gmu \|_{L^2}^2 +
		\| \phi_{\sim} \Gmu  \|_{L^2}^2 +
		\| \phi_> \Gmu  \|_{L^2}^2 .
	\]

	The first two terms are trivial to deal with.
	If $|x|\leq M^2$, we have,
	\[
		\Gmu(x) \leq \exp( \mu |x|^2 ) \leq \exp(\mu M^4),
	\]
	so setting $\mu \leq M^{-4}$ gives
	$
		\| \phi_< \Gmu \|_{L^2} \leq 
		\| \phi_< e^1 \|_{L^2} \leq e \| \phi \|_{L^2},$ uniformly in $\epsilon$.
		The same bound holds for $\phi_\sim$.
	It remains then to bound 
		$\| \phi_> \Gmu  \|_{L^2} $ uniformly in $\epsilon$.

	Starting with equation  \eqref{eqn:4:weakSW}, we multiply both sides by $\phi_>(x) \Gmu(x)^2 $,
	\[
		\omega  \phi(x) \Gmu(x)^2  \leq \T(\phi,\ldots,\phi)(x) \phi(x) \Gmu(x)^2.
	\]
	Now integrating over $\R$, using the relationship between $\Eqin$ and $\Tqin$ in \eqref{eqn:4:TintermsofE},
		and passing the exponential weight using \eqref{eqn:4:weight}, we determine the bound,
	\begin{align*}
		\omega \| \phi \Gmu \|_{L^2}^2 & \leq 6 \Eqin(\phi,\phi,\phi,\phi,\phi,\phi_>  \Gmu^2 )	
			\lesssim \Eqin( \phi e^\Gmu,\ldots,\phi e^\Gmu, \phi_> e^\Gmu).
	\end{align*}
	For convenience, let $\psi = \phi \Gmu$.
	The bound then reads,
	\begin{equation}
		\omega \| \psi \|_{L^2}^2 \lesssim \Eqin(\psi,\psi,\psi,\psi,\psi,\psi_>).
		\label{eqn:4:lastpsi}
	\end{equation}
	Now write each $\psi = \psi_< + \psi_{\sim} + \psi_>$ and expand the multinear functional.
	We will get many terms, which we bound in one of two ways.
	\begin{itemize}
		\item
			If there are three or more $\psi_>$ terms, bound by $\|\psi_>\|^k_{L^2}$ (where $k \geq 3$ is the number of $\psi_>$ terms appearing)
			using the standard $L^2$ bound \eqref{eqn:4:bound}.
			In this case the other terms are $\psi_<$ or $\psi_\sim$, which we know are uniformly bounded.
		\item
	If there are one or two $\psi_>$ terms, then there is either a $\psi_<$ term or a $\psi_\sim$ term.
			We may assume $M>4$.
			Then in the former case we can use the refined multilinear estimate  \eqref{eqn:4:refined} (with $r=M$ and $R=M^2$) 
			and bound by $(1/M^{1/3}) \|\psi_>\|_{L^2}^k$ 
		(where $k=1$ or $k=2$).
			If there are no $\psi_<$ terms, we 
			bound by $\| \psi_\sim \|_{L^2} \|\psi_> \|_{L^2}^k \lesssim \| \phi_\sim \|_{L^2} \|\psi_> \|_{L^2}^k$.
	\end{itemize}
	Using these, we find,
	\begin{align}
		\omega \| \psi_> \|_{L^2}^2
                & \leq 6 \Eqin(\psi,\psi,\psi,\psi,\psi,\psi_>)	\notag \\
		&\leq C\left(  \sum_{k=3}^m \| \psi_> \|^k + \left(\frac{1}{M^{1/3}}  + \| \phi_\sim \|_{L^2}\right) (
				\|\psi_>\|_{L^2}^2 +
				\|\psi_>\|_{L^2} ) \right), \label{eqn:4:poly}
	\end{align}
	for some constant $C$ independent of $\psi$ and $\epsilon$.
	Set $\delta(M) = M^{-1/3} + \| \phi_\sim\|_{L^2}$ 
		and 
		\[
			x(\epsilon,M) = \| \psi_>\|_{L^2} = \| \phi \chi_{|x|<M^2} \Gmu \|_{L^2}.
		\]
	Note that $\delta(M) \rightarrow 0 $ as $M \rightarrow \infty$.
	Choose $M$ sufficiently large so that 
		$C \delta(M) \leq \omega/2$.
		This gives,
	\[
		\frac{\omega}{2C} x(\epsilon,M)^2
			\leq   \sum_{k=3}^m x(\epsilon,M)^k + \delta(M) x(\epsilon,M).
	\]
	Dividing through by $x(\epsilon,M)>0$ and rearranging terms gives,
	\begin{equation}
		0 \leq p_{\delta(M)}( x(\epsilon,M)), \;\;\text{ where }\;\;
		p_\delta(x):= 
			\sum_{k=2}^{m-1} x^k -\frac{\omega}{2C} x+  \delta.
			\label{eqn:4:finalpoly}
	\end{equation}

	Observe that $p_0(0)=0$,  $p_0'(0) = -\omega/2C <0$ and  $p_0(x) \rightarrow \infty$ as $x \rightarrow \infty$.
	This shows that $p_0$ has another $0$ in $(0,\infty)$; call the smallest such zero $x_0$.
	The zeroes of a polynomial are continuous functions of the coefficients.
	Hence if we choose $M$ sufficiently large we can assume that $p_{\delta(M)}$ has one zero in $(-\infty,x_0/3)$ (coming from $p_0(0)=0)$ and 
		one zero $(2x_0/3,\infty)$ (coming from $p_0(x_0)=0$)
		and that $p_{\delta(M)}(x)<0$ in $(x_0/3,2x_0/3)$.
	This shows that for all $M$ sufficiently large,
	\[
		Z_{\delta(M)} = p_{\delta(M)}^{-1}( [0,\infty) ) \subset (-\infty,x_0/3) \cup (2x_0/3,\infty).
	\]

	Now we know from the inequality \eqref{eqn:4:finalpoly} that $x(\epsilon,M) \in Z_{\delta(M)}$ for all $\epsilon$.
	If we set $\epsilon=1$ we get,
	\[
		x(1,M) = \| \psi_> \|_{L^2} = \| \phi_> e^{\mu x^2/(1+x^2)} \|_{L^2} \leq \| \phi_>\|_{L^2} e^\mu.
	\]
	Recall that we set $\mu =  M^{-4}$, so that $\mu \lesssim 1$ and so $x(1,M) \lesssim \| \phi_> \|_{L^2}$.
	As $M \rightarrow \infty$, $\| \phi_> \|_{L^2} = \| \phi \chi_{M^2<|x|} \|_{L^2} \rightarrow 0$, 
		and hence if we take $M$ sufficiently large we will have
	$ x(1,M) \leq x_0/3$.
	But now because $x(\epsilon,M)$ depends continuously on $\epsilon$,
	and $x(\epsilon,M)\in Z_{\delta(M)}$ for all $\epsilon$, we have,
	\[
		x(\epsilon,M) = \| \psi_> \|_{L^2} \leq x_0/3,
	\]
	for all $\epsilon$.
	Taking $\epsilon \rightarrow 0$ yields $x(0,M) = 
		\| \phi_> e^{\mu x^2 } \| \leq x_0/3 <\infty,$
	which is what we wanted to prove.
\end{proof}

\begin{corr}
	Suppose that $\phi \in L^2$ is a stationary wave solution of the Hamiltonian flow associated to $\Hqin$; that is,
	$\phi$ satisfies,
	\begin{equation}
		\omega \phi (x)= \Tqin(\phi,\phi,\phi,\phi,\phi)(x),
		\label{eqn:4:strongSW}
	\end{equation}
	for some $\omega$.
	Then there exists $\alpha>0$ and $\beta>0$ such that
	$\phi e^{\alpha x^2} \in L^\infty$ and $\widehat{\phi} e^{\beta x^2} \in L^\infty$.
	As a result, $\phi$ can be extended to an entire function on the complex plane.
	\label{thm:4:smooth}
\end{corr}

\begin{proof}
	The condition \eqref{eqn:4:strongSW} implies the condition \eqref{eqn:4:weakSW} in the previous theorem,
		and hence there exists $\alpha>0$ such that $\phi e^{2\alpha x^2} \in L^2$.
		Because $\Tqin$ commutes with the Fourier transform, condition \eqref{eqn:4:strongSW} also holds
		with $\phi$ replaced by $\hat{\phi}$.
		Then, again by the previous theorem, there exists $\beta>0$ such that $\phi e^{2\beta x^2} \in L^2$.

	To turn these $L^2$ bounds into $L^\infty$ bounds, we first assume $\phi$ is Schwartz and compute,
	\begin{align*}
		\phi(x)^2 e^{ 2\alpha x^2 }	
			&= e^{ 2\alpha x^2 } \int_x^\infty \frac{d}{dt} \phi(t)^2 dt
				=	e^{2 \alpha x^2 } \int_x^\infty 2 \phi(t) \phi'(t) dt	
				\leq 2 \int_x^\infty e^{ 2 \alpha t^2 } \phi(t) \phi'(t) dt	\\
			&\leq 2 \| e^{2\alpha t^2 } \phi(t) \|_{L^2} \| \phi' \|_{L^2}	
                =2 \| e^{2\alpha t^2 } \phi(t) \|_{L^2} \| \xi \hat{\phi} \|_{L^2} 
                \leq \beta^{-1/2} \| e^{2\alpha t^2 } \phi(t) \|_{L^2} \| e^{2\beta \xi^2} \hat{\phi} \|_{L^2}	,
	\end{align*}
	which gives $\phi(x)e^{\alpha x^2} \in L^\infty$.
	Because Schwartz functions are dense in $L^2$, this holds for arbitrary $\phi \in L^2$.
	The $L^\infty$ bound for $\hat{\phi}$ follows similarly.

	Finally, using the $L^\infty$ bound $|\hat{\phi}(\xi)| \lesssim e^{-\beta \xi^2 }$,
	and the inverse Fourier transform formula,
	\[
		\phi(z) = \frac{1}{\sqrt{2\pi}} \int_\R e^{iz\xi} \hat{\phi}(\xi) d\xi.
	\]
	we can extend $\phi$ to an entire function on the complex plane.
\end{proof}

\subsection{Smoothing and further boundedness properties}

In this section we prove two further boundedness results for the operator $\Tqin$.
Unlike our previous boundedness results, which were simply inherited from the analogous results for $T_A$,
	the present results rely on additional structure in $\Tqin$.

We first strengthen items (ii) and (iii) in Theorem \ref{thm:4:Tbounds}.
Item (ii) says that $\Tqin$ is bounded from $(L^{2,\sigma})^5$ to $L^{2,\sigma}$.
Our first theorem here improves this by showing that $\Tqin$ in fact maps $(L^{2,\sigma})^5$ to $L^{2,\sigma+\delta}$
	for some $\delta>0$.
	By Fourier transform invariance, $\Tqin$ maps $(H^\sigma)^5$ to $H^{\sigma+\delta}$.

The  second result here establishes boundedness of $\Tqin$ from $(L^{\infty,s})^5$ to $L^{\infty,s}$ for any $s>1/2$.
(The analogous result for $s<1/2$ is false, by scaling.)

\begin{thm}
	For any $\sigma>0$,
	$\Tqin$ is bounded from $(L^{2,\sigma})^5$ to $L^{2,\sigma+\delta}$ with $\delta = \sigma/(1+\sigma)>0$.
	\label{thm:4:smoothing}
\end{thm}

\begin{proof}
	By duality we need to prove that for all $f_1,\ldots,f_5,g \in L^2$ with $\|f_k\|_{L^2}=\|g\|_{L^2}=1$, we have,
\[
	\left\langle
		 \Tqin\left( \langle x\rangle^{-\sigma} f_1, \ldots, \langle x \rangle^{-\sigma}f_5 \right), \langle x \rangle^{\sigma+\delta} g 
	\right\rangle_{L^2} \lesssim 1.
\]
Unpacking this using \eqref{eqn:4:TBL1}, we see that this is the same as
\begin{align}
	\int_{-1}^{1} \int_{ \R^3 } K(x,y,z,\lambda) &f_1(z-y+x) f_2(\lambda z +x ) f_3(\lambda y - y +x)  \notag\\
				& \of_4(\lambda z - y +x) \of_5(z+ \lambda y-y+x)  \og(x) dy dz dx d\lambda \lesssim 1, \label{eqn:4:nastybound}
\end{align}
where
\[
	K(x,y,z,\lambda) = \frac{ \langle x \rangle^{\sigma+\delta} }{
		\langle z-y+x \rangle^\sigma
		\langle \lambda z +x \rangle^\sigma
		\langle \lambda y -y +x  \rangle^\sigma
		\langle \lambda z - y +x \rangle^\sigma
		\langle z+ \lambda y - y +x \rangle^\sigma
	}.
\]
(The integration in $\lambda$ here is over $[-1,1]$.
	The integral in $\lambda$ for $(-\infty,-1) \cup (1,\infty)$ can be transformed into this integral by the change of variables $\lambda \mapsto 1/\lambda$.)

The overall strategy is to identify a large set on which $K$ is bounded, 
	where controlling the integral is easy,
	and 
	use a dyadic decomposition and finer bounds on $\Tqin$ to control the integral on the set where $K$ is not bounded.
On the set where $K$ is bounded the boundedness property,
\[
	\left\langle \Tqin(f_1,f_2,f_3,f_4,f_5) ,g \right\rangle \lesssim 
	\left( \prod_{k=1}^5 \|f_k\|_{L^2}^5 \right) \|g\|_{L^2},
\]
deals with the integral automatically.
So we only need to worry about the set where $K$ is unbounded.
On the unbounded piece the refined estimate,
\[
	\left\langle (\Tqin)_{|\lambda-a|<\epsilon}(f,...,f),g \right\rangle \leq \epsilon
	\left( \prod_{k=1}^5 \|f_k\|_{L^2}^5 \right) \|g\|_{L^2},
\]
will be used to gain control.
	This refined bound is a clear consequence of representation \eqref{eqn:4:EBL2}.

Fix $\epsilon$ small.
We observe first that if $\epsilon |x| \leq 1$ then $K\lesssim 1$.
So we assume that $\epsilon |x| \geq 1$.

\newcommand{\Ta}{z-y+x}
\newcommand{\Tb}{\lambda z +x}
\newcommand{\Tc}{\lambda y - y +x }
\newcommand{\Td}{\lambda z - y +x }
\newcommand{\Te}{z+ \lambda y - y +x}
Now the relation,
\[
	|z-y+x|^2+
	|\lambda z +x|^2+
	|\lambda y -y +x|^2=
	|\lambda z - y +x |^2+
	|z+ \lambda y - y +x|^2 +|x|^2,
\]
gives that,
\[
	|x| \leq \max \{ |\Ta|, |\Tb|, |\Tc| \},
\]
and hence,
\[
	K \leq \frac{ \jap{x}^\delta }{ \jap{ \Td }^\sigma \jap{\Te}^\sigma}.
\]
Now if,
\[
	|\Td| \geq \epsilon |x|	\;\text{ or }\; |\Te| \geq \epsilon |x|,
\]
then we automatically get $K \lesssim 1$ as $\delta \leq \sigma$.
Hence we assume that,
\[
	|\Td| \leq \epsilon |x|	\;\text{ and }\; |\Te| \leq \epsilon |x|.
\]


We now make some observations about how the sizes of $|y|$ and $|z|$ affect $K$.
There are four cases.

\begin{enumerate}
\item
First, suppose $|z|$ is large or comparable to $|x|$, so $|z| \geq 2 \epsilon |x|$.
Then,
\[
	2 \epsilon |x| \leq |z| \leq |\Tc| + |\Te|	\leq |\Tc| + \epsilon |x|,
\]
and hence $\epsilon |x| \leq |\Tc|$.

\item
Next, suppose $|y|$ is large or comparable to $|x|$, so $|y| \geq 2 \epsilon |x|$.
Then,
\[
	2 \epsilon |x| \leq |y| \leq |\Tb| + |\Td|	\leq |\Tc| + \epsilon |x|,
\]
and hence $\epsilon |x| \leq |\Tb|$.

\item
Next, suppose $|z|$ is small compared to $|x|$, so $|z| \leq 2 \epsilon |x|$.
Then,
\[
	|\Tb| \geq |x| - |\lambda z| \geq |x| - |z| \geq (1-2\epsilon) |x|,	
\]
and so (by the smallness of $\epsilon$) we have $\epsilon |x| \leq |\Tb|$.

\item
Finally, suppose $|y|$ is small compared to $|y|$, so $|y| \leq 2 \epsilon |x|$.
Then,
\[
	|\Tc| \geq |x| - |(\lambda-1) y| \geq |x| - 2|y| \geq (1-4\epsilon) |x|,	
\]
and so (by the smallness of $\epsilon$) we have $\epsilon |x| \leq |\Tc|$.
\end{enumerate}

From these observations we see that if $|z|$ and $|y|$ are both large we get,
\[
	\jap{x}^{\sigma+\delta} \lesssim \jap{ \Tb }^\sigma \jap{ \Tc}^\sigma,
\]
and hence $K\lesssim 1$.
If $|z|$ and $|y|$ are both small then we have the same bound on $\jap{x}$ and the same conclusion.

There are thus two regimes to consider: when $|y|$ is large and $|z|$ small, and when $|z|$ is large and $|y|$ small.
For these regimes we will use a dyadic decomposition:
\[
	\jap{x} \sim 2^j, \;\; \jap{\Td} \sim 2^k \;\text{ and }\; \jap{\Te} \sim 2^l.
\]

\emph{Regime One}: $|y|$ large, $|z|$ small.

From the observations above we have that $|\Tb| \gtrsim \epsilon 2^j$.
We have,
\[
	|\Ta| \geq |\lambda y| - |\Te| \geq \frac{\epsilon}{2} |\lambda x| \gtrsim |\lambda| 2^j,
\]
if we assume in addition that $|\lambda| 2^j \gtrsim \epsilon |\lambda| |x|/2 \geq |\Te|=2^l$.
Under this assumption we have,
\[
	K \lesssim \frac{ ( 2^j )^{\sigma+\delta} }{
		\jap{2^l} \jap{ |\lambda| 2^j } \jap{2^j} 
	} \lesssim 2^{j(\delta-\sigma)}  2^{-l\sigma} |\lambda|^{-\sigma},
\]
which is bounded if $2^{j(\delta/\sigma-1)}  2^{-l} \lesssim |\lambda|$.
Hence in Regime One, $K$ is bounded unless,
\[
	|\lambda| \leq \max\left\{ 2^{l-j},2^{j(\delta/\sigma-1)}  2^{-l}  \right\} =: \alpha.
\]

By the bound,
\[
	K \lesssim \frac{ |2^j|^{\sigma+\delta} }{ |2^k|^\sigma |2^l|^\sigma |2^j|^{\sigma} } =
		\frac{ 2^{j \delta} }{ 2^{k\sigma} 2^{l \sigma} },
\]
we then have,
\begin{align*}
	\Tqin\text{ on Regime One set} &\lesssim
	\sum_{\epsilon 2^j \geq 1}
	\sum_{\epsilon 2^j \geq 2^k}
	\sum_{\epsilon 2^j \geq 2^l}
		\frac{ 2^{j \delta} }{ 2^{k\sigma} 2^{l \sigma} }
		\langle (\Tqin)_{|\lambda|\leq \alpha}(f_1,f_2,f_3,f_4,f_5),g\rangle	\\
	&\lesssim
	\sum_{\epsilon 2^j \geq 1}
	\sum_{\epsilon 2^j \geq 2^k}
	\sum_{\epsilon 2^j \geq 2^l}
		\frac{ 2^{j \delta} }{ 2^{k\sigma} 2^{l \sigma} }
		 \max\left\{ 2^{l-j},2^{j(\delta/\sigma-1)}  2^{-l}  \right\} \lesssim 1.
\end{align*}

\emph{Regime Two}: $|z|$ large, $|y|$ small.

From the observations above we have that $|\Tc| \gtrsim \epsilon 2^j$.

We have,
\[
	|\Ta| \geq |(1-\lambda) z| - |\Td| \geq \frac{\epsilon}{2} |(1-\lambda) x| \gtrsim |1-\lambda| 2^j,
\]
if we assume in addition that $|1-\lambda| 2^j \lesssim \epsilon |1-\lambda| |x|/2 \geq |\Td|=2^k$.
Under this assumptions we have,
\[
	K \lesssim \frac{ ( 2^j )^{\sigma+\delta} }{
		\jap{2^k} \jap{ |1-\lambda| 2^j } \jap{2^j} 
	} = 2^{j(\delta-\sigma)}  2^{-k\sigma} |1-\lambda|^{-\sigma},
\]
which is bounded if $2^{j(\delta/\sigma-1)}  2^{-k} \lesssim |1-\lambda|$.
Hence in Regime Two, $K$ is bounded unless,
\[
	|1-\lambda| \leq \max\left\{ 2^{k-j},2^{j(\delta/\sigma-1)}  2^{-k}  \right\} =: \alpha.
\]

By the bound,
\[
	K \lesssim \frac{ |2^j|^{\sigma+\delta} }{ |2^k|^\sigma |2^l|^\sigma |2^j|^{\sigma} } =
		\frac{ 2^{j \delta} }{ 2^{k\sigma} 2^{l \sigma} },
\]
we then have
\begin{align*}
	\Tqin\text{ on Regime Two set} &\lesssim
	\sum_{\epsilon 2^j \geq 1}
	\sum_{\epsilon 2^j \geq 2^k}
	\sum_{\epsilon 2^j \geq 2^l}
		\frac{ 2^{j \delta} }{ 2^{k\sigma} 2^{l \sigma} }
		\langle \T_{|1-\lambda|\leq \alpha}(f,\ldots,f),g\rangle	\\
	&\lesssim
	\sum_{\epsilon 2^j \geq 1}
	\sum_{\epsilon 2^j \geq 2^k}
	\sum_{\epsilon 2^j \geq 2^l}
		\frac{ 2^{j \delta} }{ 2^{k\sigma} 2^{l \sigma} }
		 \max\left\{ 2^{k-j},2^{j(\delta/\sigma-1)}  2^{-k}  \right\}	\lesssim 1.
\end{align*}
The bound \eqref{eqn:4:nastybound} is thus established.
\end{proof}

\begin{thm}
	For all $s>1/2$, there is a constant $C$ such that,
	\begin{equation}
		\| \Tqin(f_1,f_2,f_3,f_4,f_5) \|_{L^{\infty,s}} \leq C \prod_{k=1}^5 \| f_k \|_{L^{\infty,s}}
		\label{eqn:4:Linfty}
	\end{equation}
	\label{thm:4:Linftys}
\end{thm}

\begin{lemma}
	Suppose that $v_1^2 + v_2^2 + v_3^2=1$.
	Then,
	\begin{equation}
		\jap{ v_1y_1 + v_2y_2 + v_3y_3 } \leq \sqrt{2}\left[ |v_1|\jap{y_1} + |v_2|\jap{y_2} + |v_3|\jap{y_3} \right].
		\label{eqn:4:myjap}
	\end{equation}
\end{lemma}
\begin{proof}
	We have,
	\begin{align*}
		1+ ( v_1y_1 + v_2y_2 + v_3y_3 )^2
		&\leq 2 \left(  1+ v_1^2 y_1^2 + v_2^2 y_2^2 + v_3^2 y_3^2 \right)	\\
		&= 2 \left(   v_1^2 (1+y_1^2) + v_2^2 (1+y_2^2) + v_3^2 (1+y_3^2) \right)	\\
		&\leq 2 \left(   |v_1| (1+y_1^2)^{1/2} + |v_2| (1+y_2^2)^{1/2} + |v_3| (1+y_3^2)^{1/2} \right)^{2}.
	\end{align*}
	Taking square roots then gives \eqref{eqn:4:myjap}.
\end{proof}

\begin{proof}[Proof of Theorem \ref{thm:4:Linftys}]
	We may assume by rescaling that $\| f_k \|_{L^{\infty,s}} = 1$, which means $|f_k(t)|\leq \jap{t}^{-s}$.
	Set $y_k = (A(\lambda) x)_k$.
	Then, for $x=(x_1,x_2,x_3)\in\R^3$,
	\begin{align*}
		\| \Tqin(f_1,f_2,f_3,f_4,f_5) \|_{L^{\infty,s}}
			&\leq \sup_{x_1\in \R} \left(  \jap{x_1}^s \Tqin\left( \jap{t}^{-s}, \jap{t}^{-s}, \jap{t}^{-s}, \jap{t}^{-s}, \jap{t}^{-s} \right)(x_1)\right)	\\
			&= \sup_{x_1 \in \R} \int_\R \int_{\R^2} \frac{1}{\lambda^2-\lambda+1} 
			\frac{ \jap{x_1}^s }{ \jap{y_1}^s \jap{y_2}^s \jap{y_3}^s \jap{x_2}^s \jap{x_3}^s } dx_2 dx_3 d\lambda.
	\end{align*}
	We will show that the integral of the Japanese bracket terms over $x_2$ and $x_3$ can be bounded by an absolute constant
		independent of $\lambda$.
	Because $(\lambda^2-\lambda+1)^{-1}$ is integrable over $\R$, this will prove the bound.
	
	By \CS{}, we have,
	\begin{equation}
		\begin{split}
			&\int_{\R^2} 
			\frac{ \jap{x_1}^s }{ \jap{y_1}^s \jap{y_2}^s \jap{y_3}^s \jap{x_2}^s \jap{x_3}^s } dx_2 dx_3 
			\\
			&\hspace{2cm}
                        \leq
			 \left( \int_{\R^2} \frac{ \jap{x_1}^{2s} }{ \jap{y_1}^{2s} \jap{y_2}^{2s} \jap{y_3}^{2s}  } dx_2 dx_3 \right)^{1/2}
			 \left( \int_{\R^2} \frac{ 1 }{ \jap{x_2}^{2s} \jap{x_3}^{2s} } dx_2 dx_3 \right)^{1/2}.
			 \label{eqn:4:twointegrals}
		\end{split}
	\end{equation}
	The second integral here splits as $\int_\R \jap{x_2}^{-2s} dx_2 \int_\R \jap{x_3}^{-2s} dx_3$,
		and is thus finite as $s>1/2$.

	To bound the first integral we must use some structure of $A(\lambda)$, which is given by,
	\[
		A(\lambda)=
		\frac{1}{\lambda^2-\lambda+1} 
		\begin{pmatrix}
			\lambda	&1-\lambda	&\lambda^2-\lambda	\\
			\lambda^2-\lambda	&\lambda	&1- \lambda	\\
			1-\lambda	&\lambda^2-\lambda	&\lambda
		\end{pmatrix}:=
		\begin{pmatrix}
			a_{11}	&a_{12}	&a_{13}	\\
			a_{21}	&a_{22}	&a_{23}	\\
			a_{31}	&a_{32}	&a_{33}	\\
		\end{pmatrix}.
	\]
	First we observe that the matrix has the following property.
	If we fix a row $k$ and column $j$, then the determinent of the matrix 
		obtained by deleting row $k$ and column $j$ is precisely $a_{kj}$ -- the element in row $k$ and column $j$.
	This means that,
	\begin{equation}
		a_{11} = \frac{ \lambda }{ \lambda^2 - \lambda + 1 } = \det\left[ \frac{1}{\lambda^2 - \lambda + 1 } 
			\begin{pmatrix} \lambda	&1-\lambda \\ \lambda^2-\lambda	& \lambda\end{pmatrix} \right]
				= \det \begin{pmatrix}  a_{22} & a_{23} \\ a_{32} & a_{33} \end{pmatrix},
					\label{eqn:4:amazingdet}
	\end{equation}
	with similar formulas for $a_{21}$ and $a_{31}$.

	Next, because $A$ is an isometry, the inverse matrix is just the transpose.
	Since $x = A^{-1}y$ we have the formula,
	$
		x_1 
			= a_{11} y_1 + a_{21} y_2 + a_{31} y_3.
	$
	$A$ is an isometry, so $a_{11}^2 + a_{21}^2 + a_{31}^2=1$.
	We can therefore use \eqref{eqn:4:myjap}, from the previous lemma, raised to to the power $s$; it reads,
	$
		\jap{x_1}^{2s} \lesssim
			|a_{11}|^{2s} \jap{y_1}^{2s} +
			|a_{21}|^{2s} \jap{y_2}^{2s} +
			|a_{31}|^{2s} \jap{y_3}^{2s}.
	$
	Applying this to bound the first integral in \eqref{eqn:4:twointegrals}, we then have,
	\begin{align*}
		&\int_{\R^2} \frac{ \jap{x_1}^{2s} }{ \jap{y_1}^{2s} \jap{y_2}^{2s} \jap{y_3}^{2s}  } dx_2 dx_3 	
                \\
		 &\hspace{2cm}
                 \lesssim
		 \int_{\R^2} \frac{ |a_{11}| }{  \jap{y_2}^{2s} \jap{y_3}^{2s}  } dx_2 dx_3 +
		 \int_{\R^2} \frac{ |a_{21}| }{  \jap{y_1}^{2s} \jap{y_3}^{2s}  } dx_2 dx_3 + 
		 \int_{\R^2} \frac{ |a_{31}| }{  \jap{y_1}^{2s} \jap{y_2}^{2s}  } dx_2 dx_3.
	\end{align*}
	We will show how the first integral may be bounded; the other two are bounded by an identical argument.
	We perform the change of variables $z_2 = y_2 = (Ax)_2$ and $z_3 = y_3 = (Ax)_3$.
	Expressed as a matrix, this change of variables is,
	\[
		\begin{pmatrix} z_2 \\ z_3 \end{pmatrix}
		=\begin{pmatrix}  a_{22} & a_{23} \\ a_{32} & a_{33} \end{pmatrix}
		\begin{pmatrix} x_2 \\ x_3 \end{pmatrix}.
	\]
	The determinent of this change of variables is, by \eqref{eqn:4:amazingdet}, simply $|a_{11}|$.
	Therefore, using that $|a_{11}|\leq1$ because $A$ is an isometry,
	\[
		\int_{\R^2} \frac{ |a_{11}|^{2s} }{  \jap{y_2}^{2s} \jap{y_3}^{2s}  } dx_2 dx_3 =
		 \int_{\R^2} \frac{ |a_{11}|^{2s-1} }{  \jap{z_2}^{2s} \jap{z_3}^{2s}  } dz_2 dz_3 \leq
		 \int_{\R} \frac{ 1 }{  \jap{z_2}^{2s}   } dz_2 
		 \int_{\R} \frac{ 1 }{  \jap{z_3}^{2s}   } dz_3 ,
	\]
	and the right hand side is finite because $s>1/2$.	
\end{proof}

\section{The cubic resonant equation}

In this final section we study the system defined by the Hamiltonian,
\begin{equation}
		\Hcub(f) = \frac{2}{\pi} \| e^{itH} f \|_{L^4_t L^4_x  }^4
	= \frac{2}{\pi} \int_0^{\pi/2} \int_\R | e^{itH}f(x) |^4 dx dt,
	\label{eqn:5:Hhermite}
\end{equation}
which has an associated multilinear functional,
\begin{equation}
	\Ecub(f_1,f_2,f_3,f_4) 
	= \frac{2}{\pi} \int_0^{\pi/2} \int_\R  
		 (e^{itH}f_1(x)) ( e^{itH}f_{2}(x) )
		 \overline{(e^{itH}f_3(x)) (e^{itH}f_{4}(x) }) dx dt .
	\label{eqn:5:Ehermite}
\end{equation}
The functional has a large number of permutation symmetries,
\begin{equation}
	\Ecub(f_1,f_2,f_3,f_4)=
	\Ecub(f_2,f_1,f_3,f_4)=
	\Ecub(f_1,f_2,f_4,f_3)=
	\overline{\Ecub(f_3,f_4,f_1,f_2)}.
\end{equation}
As in the case of quintic resonant system, Hamilton's equation is $iu_t = \Tcub(u,u,u)$ where $\Tcub$ is defined by,
\begin{equation}
	\jap{ \Tcub(f_1,f_2,f_3),g } = 4 \Ecub(f_1,f_2,f_3,g).
	\label{eqn:5:TintermsofE}
\end{equation}
By an identical computation to the derivation of \eqref{eqn:4:hameqn},
we find that the operator $\Tcub$ is given explictely by,
\begin{equation}
	\Tcub(f_1,f_2,f_3)(x)
		= \frac{8}{\pi} \intpi 
			e^{-itH}\left[ (e^{itH}f_1)(e^{itH}f_2))
			\overline{(e^{itH}f_3)} \right] (x) dt.
	\label{eqn:5:Thermite}
\end{equation}
This shows that the flow corresponding to the Hamiltonian $\Hcub$ is the resonant equation \eqref{eqn:2:respde2} in the case $k=1$, up to rescaling time.
As discussed in the introduction, it was shown in \cite{Hani2015} that this system is also the modified scattering limit
	of the NLS equation \eqref{eqn:1:hani}.

\subsection{Representations of the Hamiltonian and the flow operator}

As for the quintic case, we devote a significant amount of work to determining alternative representations of $\Ecub$, $\Hcub$ and $\Tcub$.
In contrast to the quintic case, we do \emph{not} have representations for $\Hcub$ of the form,
\begin{align*}
	\int_{\R^4} f_1(y_1) f_2(y_2) \overline{ f_3(y_3) f_4( y_4) } \delta_{y_1+y_2=y_3+y_4} \delta_{ y_1^2+y_2^2=y_3^2+y_4^2} dy
	\;\;\text{ or }\;\;\| e^{it\Delta}f \|_{L^4_t L^4_x}^4.
\end{align*}
These representations are inconsistent with the scaling of the inequality 
	$\Hcub(f) \leq (1/\sqrt{8\pi}) \| f \|_{L^2}^4$
	which we prove in Theorem \ref{thm:5:L2}.

\begin{thm}
	There holds the representations,
	\begin{align}
		\Ecub(f_1,f_2,f_3,f_4)
		&=
		\frac{1}{2\pi^2}
		\int_{\R^4}
		 e^{-\frac{1}{2}[ (\lambda v_2 )^2 + v_1^2 ] }
		f_1(\lambda v_1 + v_3) f_2(v_2+v_3) \notag \\
		&\hspace{5cm} \overline{ f_3(\lambda v_1 + v_2 + v_3 )f_4(v_3) }
		dv_1 dv_2 dv_3 d\lambda,
		\label{eqn:5:EBL1}		\\
		\Tcub(f_1,f_2,f_3)(x)
		&=
		\frac{2}{\pi^2}
		\int_{\R^4}
		 e^{-\frac{1}{2}[ (\lambda v_2 )^2 + v_1^2 ] }
		f_1(\lambda v_1 + x) f_2(v_2+x)
		\overline{ f_3(\lambda v_1 + v_2 + x ) }
		dv_1 dv_2 d \lambda,
		\label{eqn:5:TBL1}
	\end{align}
		\label{thm:5:BL1}
\end{thm}

(Compare with Theorem \ref{thm:4:EBL1}.)
To prove this theorem we need a lemma.

\begin{lemma}
	Let $\psi(x) = (1+x^2)^{-1/2}$. 
	Then $\hat{\psi}(\xi) = \zeta(\xi) := \displaystyle \frac{1}{\sqrt{2\pi}} \int_\R \frac{1}{|v|} e^{-\frac{1}{2}[ (\xi/v)^2 + v^2 ] } dv$.
	\label{thm:5:fourierlemma}
\end{lemma}

\begin{proof}
	It is clear that $\psi$ is in $L^2$.
	We will calculate the Fourier transform of  $\zeta(\xi)$
	and find that it equals $\psi$.
	The lemma then follows from Fourier inversion and the fact that $\psi$ is even.

	We have,
	\begin{align*}
		\hat{\zeta}(x) 
			&= 
		    \frac{1}{ {2\pi} } \int_\R e^{ix\xi}
	    	\int_\R \frac{1}{|v|} e^{-\frac{1}{2}[ (\xi/v)^2 + v^2 ] } dv d\xi
                \\&
			=
			\frac{1}{ {2\pi} } \int_\R \frac{1}{|v|} e^{-\frac{1}{2} v^2  } 
			\int_\R e^{ix\xi}
			e^{-\frac{1}{2} (\xi/v)^2}   d\xi dv	
			=
			\frac{1}{\sqrt{2\pi}} \int_\R e^{-\frac{1}{2} v^2  } e^{-\frac{1}{2}v^2 x^2} dv,
	\end{align*}
	where in the last equality we used the explicit Fourier transform of the Gaussian $e^{-ax^2}$ with $a=1/(2v^2)$.
	In this last integral we perform the change of variables $u = v(1+x^2)^{-1/2}$, 
        which gives the result.
	\end{proof}

	\begin{proof}[Proof of Theorem \ref{thm:5:BL1}]
		We evaluate \eqref{eqn:5:Ehermite}
                using the Mehler formula \eqref{eqn:2:mehler}, which reads,
		\[
			e^{itH}f_k(x) = \frac{1}{\sqrt{2\pi}|\sin(2t)|^{1/2}}
			\int_\R e^{ - i \frac{ (x^2/2+y^2/2)\cos(2t) -xy }{\sin(2t)} }f_k(y) dy.
		\]
		For notational convenience, let $\Lambda(x,t)= (e^{itH}f_1)(e^{itH}f_2)\overline{(e^{itH}f_3)(e^{itH}f_4)}$ be the integrand in \eqref{eqn:5:Ehermite}.
		Using the Mehler formula, we have,
\[
	\Lambda(x,t) = 
	 \frac{1}{4\pi^2 |\sin(2t)|^{2}}
		\int_{\R^4} e^{-i \frac{ \Omega \cos(2t) }{2\sin(2t)} }
		e^{-i\frac{(y_1+y_2-y_3-y_4)x}{\sin(2t)}} f_1(y_1)f_2(y_2) \overline{f_3(y_3) f_4(y_4)} dy_1dy_2dy_3dy_4,
\]
where $\Omega = y_1^2+y_2^2-y_3^2-y_4^2$.
Changing variables $w(y_3) = -y_1- y_2 + y_3 + y_4$ and integrating over $x$ yields,
\begin{align*}
	\int_\R 
		\Lambda(x,t)dx	
		&= 
	 \frac{1}{4\pi^2|\sin(2t)|^{2}} \int_\R
		\int_{\R^4} e^{-i \frac{ \Omega \cos(2t) }{2\sin(2t)} }
		e^{i\frac{wx}{\sin(2t)}} f_1(y_1)f_2(y_2) \\
                &\hspace{5cm}\overline{f_3(w + y_1 + y_2 - y_4) f_4(y_4)} dw dy_1 dy_2 dy_4 dx	\\
	 &=\frac{1}{2\pi|\sin(2t)|}
		\int_{\R^3} e^{-i \frac{ \Omega \cos(2t) }{2\sin(2t)} }
		f_1(y_1)f_2(y_2) \overline{f_3(y_1+y_2-y_4) f_4(y_4)}  dy_2dy_3dy_4,
\end{align*}
		where to get the second equality we used the Fourier inversion formula \eqref{eqn:1:fourierinv} 
with $a=1/\sin(2t)$.

We now integrate $t$ on the interval $[-\pi/4,\pi/4]$ and change of variables $u = -\cos(2t)/\sin(2t)$.
		This change of variables bijectively maps $(-\pi/4,0) \cup (0,\pi/4]$ to $(-\infty,+\infty)$ and satisfies $du = 2 dt/\sin^2(2t)$.
Moreover,
$
		u^2 =\cos^2(2t) / \sin^2(2t) = (1/\sin^2(2t)) - 1,
$
which gives $\sin(2t) = (u^2+1)^{-1/2}$.
Using these, we find,
		\begin{align*}
	\intpi
	\int_\R 
			&\Lambda(x,t) dx dt
                         = 
	\frac{1}{4\pi}
		\int_{\R^3}
		\left( \int_{\R} e^{-i \frac{ \Omega  }{2}u}
		\frac{1}{(1+u^2)^{1/2}} du \right)
		f_1(y_1)f_2(y_2) \overline{f_3(y_1+y_2-y_4) f_4(y_4)} 
                dy_2dy_3dy_4.\\
                &=
	\frac{1}{4\pi}
		\int_{\R^3}
		\left( \int_{\R} \frac{1}{|v_1|} e^{-\frac{1}{2}[ (\Omega/2v_1)^2 + v_1^2 ] }
		dv_1
		\right)
		f_1(y_1)f_2(y_2) \overline{f_3(y_1+y_2-y_4) f_4(y_4)} 
                dy_2dy_3dy_4,
		\end{align*}
                where in the second inequality we used
                Lemma  \ref{thm:5:fourierlemma}.

At this point $\Omega/2 = [(y_1)^2 + (y_2)^2 - (y_1+y_2-y_4)^2 - (y_4)^2 ]/2=(y_1-y_4)(y_2-y_4)$.
For fixed $v_1$, we perform the linear change of variables,
\[
	\begin{pmatrix}
		\lambda \\ v_2 \\ v_3
	\end{pmatrix} =
	\begin{pmatrix}
		(y_1-y_4)/v_1 \\ y_2-y_4 \\ y_4
	\end{pmatrix} =
	\begin{pmatrix}
		1/v_1	&0	&- 1/v_1	\\
		0	&1	&-1	\\
		0	&0	&1
	\end{pmatrix}
	\begin{pmatrix}
		y_1 \\ y_2 \\ y_4
	\end{pmatrix},
\]
which has determinent $1/|v_1|$. 
The inverse is given by,
\[
	\begin{pmatrix}
		y_1 \\ y_2 \\ y_4
	\end{pmatrix} = 
	\begin{pmatrix}
		\lambda v_1 + v_3  \\ v_2 + v_3	\\
		v_3 
	\end{pmatrix} =
	\begin{pmatrix}
		v_1	&0	&1	\\
		0	&1	&1	\\
		0	&0	&1	
	\end{pmatrix}
	\begin{pmatrix}
		\lambda \\ v_2 \\ v_3
	\end{pmatrix} ,
\]
and we note specifically that $\Omega/2v_1 = [(y_1-y_4)/v_1](y_2-y_4) = \lambda v_2$.
Performing this change of variables gives \eqref{eqn:5:EBL1}.
To get \eqref{eqn:5:TBL1}, we simply use the relation $\langle \Tcub(f_1,f_2,f_3), g \rangle = 4 \Ecub(f_1,f_2,f_3,f_4)$.
\end{proof}

\begin{thm}
	Let $G(x) = e^{-\frac{1}{2}x^2}$. 
	There holds the representations,
	\begin{align}
		\Ecub(f_1,f_2,f_3,f_4) &= \frac{1}{\sqrt{2} \pi^2 } \int_\R \frac{1}{1+\lambda^2} E_{B(\lambda)}(G,f_1,f_2,G,f_3,f_4) d\lambda,	\label{eqn:5:EBL2}	\\
		\Tcub(f_1,f_2,f_3) &= \frac{1}{\sqrt{2} \pi^2 } \int_\R \frac{1}{1+\lambda^2} T_{B(\lambda)}(G,f_1,f_2,G,f_3) d\lambda.	\label{eqn:5:TBL2}	
	\end{align}
	where for every $\lambda$, $B(\lambda)$ is an isometry  and  $B(\lambda)(0,1,1)=(0,1,1)$.
\end{thm}

\begin{proof}
	We observe that,
	$
		(\lambda v_2)^2 + (v_1)^2 =
			\left( ( \lambda v_2 - v_1 )/ \sqrt{2}  \right)^2 
			+ \left( ( \lambda v_2 + v_1 )/ \sqrt{2}  \right)^2 
	$
	which gives,
	\[
		e^{ -\frac{1}{2} [ (\lambda v_2)^2 + (v_1)^2 ] } = 
			G\left( \frac{ \lambda v_2 - v_1 }{ \sqrt{2} } \right)
			G\left( \frac{ \lambda v_2 + v_1 }{ \sqrt{2} } \right).
	\]
	We substitute this expression into \eqref{eqn:5:EBL1}.
	Using the fact that $\overline{G(x)}=G(x)$, this gives,
	\begin{align*}
		\Ecub(f_1,f_2,f_3,f_4)	
		&=  \frac{1}{2\pi^2}
		\int_{\R^4}
			G\left( \frac{ \lambda v_2 + v_1 }{ \sqrt{2} } \right)
		f(\lambda v_1  + v_3 )f(v_2 + v_3)\\
		&\hspace{3cm}
                \overline{ G\left( \frac{ \lambda v_2 - v_1 }{ \sqrt{2} } \right) }
		\overline{f(\lambda v_1 +v_2+ v_3) f(v_3)}  dv_1 dv_2 dv_3.
	\end{align*}
	By looking at the arguments of the functions in the integrand, 
		we are led to define the matrices $C(\lambda)$ and $D(\lambda)$ by,
	\[
		C(\lambda) 
		\begin{pmatrix}
			v_1 \\ v_2 \\ v_3
		\end{pmatrix} 
		=
		\begin{pmatrix}
			1/\sqrt{2}	&\lambda/\sqrt{2}	&0	\\
			\lambda		&0	&1	\\
			0	&1	&1
		\end{pmatrix} 
		\begin{pmatrix}
			v_1 \\ v_2 \\ v_3
		\end{pmatrix} 
                \;\;\text{ and }\;\;
		D(\lambda) 
		\begin{pmatrix}
			v_1 \\ v_2 \\ v_3
		\end{pmatrix} 
		=
		\begin{pmatrix}
			-1/\sqrt{2}	&\lambda/\sqrt{2}	&0	\\
			\lambda		&1	&1	\\
			0	&0	&1
		\end{pmatrix} 
		\begin{pmatrix}
			v_1 \\ v_2 \\ v_3
		\end{pmatrix}. 
	\]
	We perform the change of variables $w = D(\lambda) v$,
        and set $B(\lambda) = C(\lambda)D(\lambda)^{-1}$.
        An identical process to the proof of Theorem \eqref{thm:4:EBL1} then
        gives \eqref{eqn:5:EBL2}.
	A calculation reveals that $B(\lambda)$ is given explicitly by,
	\begin{equation}
		B(\lambda) = C(\lambda)D(\lambda)^{-1} =
		\frac{1}{1 + \lambda^2}
		\begin{pmatrix}
			-1 + \lambda^2		&\lambda\sqrt{2} &-\lambda\sqrt{2}	\\
			-\lambda\sqrt{2}	&\lambda^2	&1	\\
			\lambda\sqrt{2}		&1	&\lambda^2 
		\end{pmatrix}.
		\label{eqn:5:Blambda}
	\end{equation}
\end{proof}

\begin{thm}
	There holds the representations,
	\begin{align}
		\E_4(f_1,f_2,f_3,f_4) &=  \frac{ 1 }{ 2\sqrt{2}  \pi^2} \int_{0}^{2\pi} E_{S(\theta)}(G,f_1,f_2,G,f_3,f_4) d\theta,	
		\label{eqn:5:EBL3}		\\
		\T_4(f_1,f_2,f_3)(x) &= \frac{\sqrt{2}}{\pi^2} \int_0^{2\pi} T_{S(\theta)}(G,f_1,f_2,G,f_3)(x) d\theta,
		\label{eqn:5:TBL3}
	\end{align}
	where $S(\theta)$ is the rotation of $\R^3$ by $\theta$ radians about the axis $(0,1,1)$.
\end{thm}

\begin{proof}
	Because the matrix $B(\lambda)$ is an isometry, $\det(B(\lambda))=+1$, and $B(\lambda)(0,1,1) = (0,1,1)$,
	it must, in fact, be a rotation about the axis $(0,1,1)$.
        An identical process to the proof of Theorem \ref{thm:4:EBL3}
        then gives the formulae.
\end{proof}


\subsection{Symmetries of the Hamiltonian and conserved quantities of the flow}

\begin{thm}
	\label{thm:5:symE}
	The function $\Ecub(f_1,f_2,f_3,f_4)$ is invariant under the following actions:
	\begin{enumerate}[label=(\roman*)]
		\item	Fourier transform, $f_k \mapsto \ft_k$.
		\item	Modulation, $f_k \mapsto e^{i\lambda} f_k$.
		\item	Linear modulation, $f_k \mapsto e^{i\lambda} f_k$.
		\item	Translation, $f_k \mapsto f_k(\cdot+\lambda)$.
		\item	\Sch{} with harmonic trapping group, $f_k \mapsto e^{i\lambda H} f_k$.
	\end{enumerate}
\end{thm}

\begin{proof}
	Because $S(\theta)$ is an isometry for all $\theta$,
        the properties of $E_A$ determined in the proof of Theorem \ref{thm:4:symE} 
        apply here too.
        Because $G$ is invariant under the Fourier transform 
        and the action $G \mapsto e^{i\lambda H}G = e^{it}G$,
        the symmetries of $E_{S(\theta)}$ carry 
        are carried over to $\Ecub$.
\end{proof}

\begin{corr}
	We have the following commuter equalities,
	\begin{align}
		e^{i\lambda Q} \Tcub(f_1,f_2,f_3) &= \Tcub(e^{i \lambda Q}f_1,e^{i\lambda Q}f_2,e^{i\lambda Q}f_3) 	
			\label{eqn:5:Tcomm1} \\
		Q \Tcub(f_1,f_2,f_3) &= \Tcub(Qf_1,f_2,f_3) + 
			\Tcub(f_1,Qf_2,f_3) - \T(f_1,f_2,Qf_3). 		\label{eqn:5:Tcomm2} 
	\end{align}
	where $Q$ are the operators: $Q=1$, $Q=x$, $Q=id/dx$, and $Q=H$.
\end{corr}

The Corollary follows immediately from Theorem \ref{thm:5:symE}
in the same way as Corollary \ref{thm:4:TAcomm1}.
By Noether's Theorem,
we determine four  conserved quantities for the Hamiltonian flow corresponding to $\Hcub$.
These are summarized in Table \ref{tbl:5:conserved}.

\begin{table}
        \caption{Symmetries of $\Hcub$ and conserved quantities of the 
				\label{tbl:5:conserved}
        cubic resonant equation.}
\begin{center}
	\renewcommand{\arraystretch}{1.5}
	\begin{tabular}{|ccc|}\hline
		Symmetry of $\Hcub$		&Conserved quantity	&Operator commuting with $\Tcub$\\	\hline
				$f\mapsto e^{i\lambda}f$	&$\int_\R |f(x)|^2 dx$		&1	\\
				$f\mapsto e^{i\lambda x}f$	&$\int_\R x|f(x)|^2 dx$	&$x$	\\
				$f\mapsto f(\cdot+\lambda)$	&$\Re\int_\R f'(x) \of(x)   dx$	 &$d/dx$\\
				$f\mapsto e^{i\lambda \HO}f$	&$\int_\R |x f(x)|^2 + |f'(x)|^2 dx$ &$H$	\\ \hline
	\end{tabular}
\end{center}
\end{table}

\subsection{Boundedness of the functional and wellposedness of Hamilton's equation}

\begin{prop}
	We have the following sharp bound,
	\begin{equation}
		|\Ecub(f_1,f_2,f_3,f_4)| \leq \frac{1}{\sqrt{2\pi}} 
			\prod_{k=1}^4 \| f_k \|_{L^2},
			\label{eqn:5:bound}
	\end{equation}
	with equality if and only if the functions $f_k$ are the same Gaussian $f_k(x) = e^{-\frac{1}{2}x^2+\beta x}$ for some $\beta \in \C$.

	In particular there holds $\Hcub(f) \leq (1/\sqrt{2\pi}) \|f\|_{L^2}^4$
	with equality if and only if $f(x) = e^{-\frac{1}{2}x^2 + \beta x}$ for some $\beta \in \C$.
	\label{thm:5:L2}
\end{prop}

The equality case here is a little different to the analogous result for $\Eqin$ in Theorem \ref{thm:4:L2bound}.
For $\Eqin$, the set of saturating functions is all Gaussians of the form $e^{-\alpha x^2 + \beta x}$ with $\Re \alpha >0$.
In the case of $\Ecub$, we necessarily have $\alpha=1/2$.

\begin{proof}
	Using the representation \eqref{eqn:5:EBL3}, we find that,
	\begin{align*}
		|\Ecub(f_1,f_2,f_3,f_4)|
		&\leq  \frac{1}{ 2 \sqrt{2} \pi^2 }
		\int_0^{2\pi} |E_{S(\theta)} (G,f_1,f_2,G,f_3,f_4)| d\theta,	\leq  \frac{1}{  \sqrt{2} \pi }
		  \| G \|_{L^2}^2 
			\prod_{k=1}^4 \| f_k \|_{L^2} ,
	\end{align*}
	We calculate
	$
		\| G \|_{L^2}^2 = \int_\R ( e^{-x^2/2} )^2 dx = \sqrt{\pi},
	$
	which yields the inequality.

    The analysis of the equality case is similar to that of Theorem \ref{thm:4:L2bound}.
    The only difference is that while in Theorem \ref{thm:4:L2bound}, the condition
    was that $f_1,\ldots,f_6$ must be the same Gaussian, here
    $f_1,\ldots,f_4$ must be the same Gaussian and equal to $G$
    up to linear modulation.
    This accounts for the restriction that $\alpha = 1/2$ in the 
    saturating Gaussian.
\end{proof}

\begin{thm}
	We have the operator bound
	$
		\| \Tcub(f_1,f_2,f_3) \|_X \leq C_X \prod_{k=1}^3 \| f_k \|_X,
	$
	for the spaces:
	\begin{enumerate}[label=(\roman*)]
		\item $X=L^2$ with $C_X=\sqrt{8/\pi}$.
		\item $X=L^{2,\sigma}$, for any $\sigma \geq 0$.
		\item $X=H^{\sigma}$, for any $\sigma \geq 0$.
		\item $X=L^{\infty,s}$, for any $s > 1/2$.
		\item $X=L^{p,s}$, for any $p \geq 2$ and $s> 1/2 - 1/p$.
	\end{enumerate}
\end{thm}

\begin{proof}
	The bounds (i) through (iii) follow as in
        the proof of Theorem \ref{thm:4:Tbounds},
		noting that in all cases $\| G \|_X < \infty$.
	
	For (iv), we need to  show $\sup_{x \in \R} |\Tcub( \jap{t}^{-s},\jap{t}^{-s},  \jap{t}^{-s})(x) \jap{x}^{s}|<\infty$.
	As in the proof of Theorem \ref{thm:4:Linftys}, it is sufficient to show that,
	\[
		\sup_{x \in \R} T_{B(\lambda)}( e^{-t^2/2}, \jap{t}^{-s}, \jap{t}^{-s}, e^{-t^2/2}, \jap{t}^{s})(x) \jap{x}^{-s} \leq C,
	\]
	for some $C$ independent of $\lambda$.
	We observe that we have $e^{-t^2/2} \lesssim \jap{t}^{-s}$,
		which means it is sufficient to show that,
	\[
		\sup_{x \in \R} T_{B(\lambda)}( \jap{t}^{-s}, \jap{t}^{-s}, \jap{t}^{-s}, \jap{t}^{-s} , \jap{t}^{-s})(x) \jap{x}^{s} \leq C,
	\]
	for some $C$ independent of $\lambda$.
	The proof of this bound is identical to the proof of the analogous bound in Theorem \ref{thm:4:Linftys}. 

	Item (v) follows from interpolating between $L^{2,\sigma}$ and $L^{\infty,s}$.
\end{proof}

\begin{thm}
	\label{thm:5:wellposedness}
	Consider the Cauchy problem,
	\begin{equation}
		\begin{split}
			iu_t &= \Tcub(u,u,u),		\\
			f(t=0) &= f_0,
		\end{split}
		\label{eqn:5:cauchy2}
	\end{equation}
	which is Hamilton's equation corresponding to $\Hcub$
	and the resonant equation \eqref{eqn:1:respde2} in the cubic case $k=1$.
	
	\begin{enumerate}[label=(\roman*)]
		\item The Cauchy problem \eqref{eqn:5:cauchy2} is locally wellposed in $X$ for any of the spaces in the previous theorem.
		\item The Cauchy problem \eqref{eqn:5:cauchy2}  is globally wellposed in $L^2$
                \item Persistance of regularity:
                for $H^\sigma$ initial data the $L^2$ global solution is in $H^\sigma$
                for all time.
	\end{enumerate}
	
\end{thm}

	The proof is identical to that of Theorem \ref{thm:4:wellposedness}.

\subsection{Analysis of the stationary waves}

Stationary waves $\psi$ for the cubic equation are solutions of the equation
    $\omega \psi = \Tcub(\psi, \psi, \psi)$ for some $\omega \in \R$.
As for the quintic equation, one may show that $\Ecub(\phi_{n_1},\phi_{n_2},
\phi_{n_3},\phi_{n_4}) = 0$ unless $n_1+n_2=n_3+n_4$.
From the definition of $\Tcub$ by duality in \eqref{eqn:5:TintermsofE},
we then see that $\phi_n$ is a stationary wave of the cubic
resonant equation for all $n \geq 0$.

By applying the symmetries of $\Hcub$, we find that all functions of the form,
\begin{equation}
	a e^{ibx} \phi_n(x+c),
	\label{eqn:5:allstats}
\end{equation}
are stationary waves for $a\in \C$ and $b,c\in \R$.
The set of stationary waves we can construct for the cubic
case is smaller than the set we can construct
for the quintic case in \eqref{eqn:4:allstats},
	because the cubic equation has fewer symmetries.

\subsubsection{Regularity of stationary waves: technical issues}

All of the stationary waves constructed in the previous subsection are analytic and exponentially decaying in space.
In the remainder of this section we prove that all stationary waves that are in $L^2$ are automatically analytic and decay in space like $e^{- \alpha x^2}$ for some $\alpha>0$.
This is analogous to Corollary \ref{thm:4:smooth} for the quintic resonant equation.
Recall that the proof of that result relied on two ingredients:
	a refined multilinear Strichartz estimate \eqref{eqn:4:refined} and a weight transfer property \eqref{eqn:4:weight}.

The weight transfer property for the quintic equation used
the analogous property for the functionals $E_A$ given in \eqref{eqn:4:EAweight}.
In the present case we encounter a problem when trying to replicate this: when we try to transfer weight in the functional $\Ecub$ in the same way,
	the weight also hits the Gaussians,
\begin{equation}
	E_{B(\lambda)}( G,f_1,f_2,G,f_3,f_4 \Gmu ) \leq
	E_{B(\lambda)}( G\Gmu,f_1\Gmu ,f_2\Gmu ,G\Gmu ,f_3\Gmu ,f_4 ),
	\label{eqn:5:weight1}
\end{equation}
and the right hand side here can't be related back to $\Ecub$.
To get around this, we observe that,
\begin{equation}
	G(x) \Gmu(x) = e^{-\frac{1}{2}x^2} e^{\mu x^2/(1+\epsilon x^2)} \leq e^{ \left( -\frac{1}{2} + \mu \right)x^2},
	\label{eqn:5:weight2}
\end{equation}
which, if $\mu<1/2$, is still decaying exponentially fast, and should be possible to handle in estimates.

Because of this consideration, we are led to define,
\begin{equation}
	\Ecub^\mu(f_1,f_2,f_3,f_4) = \frac{1}{\sqrt{2}\pi^2} \int_\R \frac{1}{1+\lambda^2}
		E_{B(\lambda)}\left(
			e^{\left(-\frac{1}{2}+\mu\right)x^2}, f_1,f_2,
			e^{\left(-\frac{1}{2}+\mu\right)x^2}, f_3, f_4 \right) d\lambda,
	\label{eqn:5:weight3}
\end{equation}
and we note that $\Ecub^0 = \Ecub$.
We now proceed to develop the two ingredients
	for the stationary wave result, noting that both ingredients need to be developed for $\Ecub^\mu$ and not just $\Ecub$.

\subsubsection{The weight transfer property}

\begin{lemma}

	(i) If $\mu<\frac{1}{2}$ and functions $f_k$ are positive then there holds,
	\begin{equation}
		\Ecub \left(f_1,f_2,f_3,f_4 \Gmu \right)
			\leq \E^\mu \left(f_1 \Gmu(x), f_2 \Gmu(x), f_3 \Gmu(x), f_4  \right).
			\label{eqn:5:Emubound1}
	\end{equation}

	(ii) If $\mu<\frac{1}{2}$ then there holds the bound,
	\begin{equation}
		|\Ecub^\mu(f_1,f_2,f_3,f_4)| \leq \sqrt{ \frac{\pi}{8} } \frac{1}{\sqrt{1-2\mu} } 
			\|f_1\|_{L^2} \|f_2\|_{L^2} \|f_3\|_{L^2} \|f_4\|_{L^2}.
			\label{eqn:5:Emubound2}
	\end{equation}
\end{lemma}


\begin{proof}
	(i) This is immediate from the computations in \eqref{eqn:5:weight1} and \eqref{eqn:5:weight2}.

	(ii)
	Boundedness is proved in the usual way,
	\begin{align*}
		|\Ecub^\mu \left(f_1,f_2,f_3,f_4 \right) |
			&\leq\frac{1}{\sqrt{8}\pi} \int_\R \frac{1}{1+\lambda^2} \left\| e^{\left( -\frac{1}{2} + \mu \right)x^2 } \right\|_{L^2}^2 
			\|f_1\|_{L^2} \|f_2\|_{L^2} \|f_3\|_{L^2} \|f_4\|_{L^2} d\lambda	\\
			&\leq\frac{1}{\sqrt{8}\pi} \left( \int_\R \frac{1}{1+\lambda^2} d\lambda \right) \left( \int_\R e^{-(1-2\mu)x^2} dx  \right)
			\|f_1\|_{L^2} \|f_2\|_{L^2} \|f_3\|_{L^2} \|f_4\|_{L^2} .
	\end{align*}
	Evaluating the integrals appearing here yields the result.
\end{proof}

\subsubsection{Refined multilinear estimates}

As for the quintic resonant equation, the refined multilinear estimates we need can be determined in an elementary way
	using the representations \eqref{eqn:5:EBL1} and \eqref{eqn:5:EBL2} for $\Ecub$.

\begin{lemma}
	There is an absolute constant $C$ such that if 
		$f_1$ and $f_3$ are supported in $B(0,R)^C$ and $f_2$ and $f_4$ are supported in $B(0,r)$, with $R>4r$, then,
	\begin{equation}
		| \Ecub(f_1,f_2,f_3,f_4) |\leq 
			\frac{C}{ \sqrt{R} }
			\| f_1 \|_{L^2}
			\| f_2 \|_{L^2}
			\| f_3 \|_{L^2}
			\| f_4 \|_{L^2}
			\label{eqn:5:refined1}
	\end{equation}
\end{lemma}

\begin{proof}
	From \eqref{eqn:5:EBL1} we have,
	\begin{align}
		\Ecub(f_1,f_2,f_3,f_4) 
			&=\frac{1}{2\pi^2}
				\int_{\R^4}
				 e^{ - \frac{1}{2} [ (\lambda v_2 )^2 + v_1^2 ] }
				f_1(\lambda v_1 + v_3) f_2(v_2+v_3) \notag\\
                            &\hspace{3cm}
				\overline{ f_3(\lambda v_1 + v_2 + v_3 )f_4(v_3) }
				dv_1 dv_2 dv_3 d\lambda. \label{eqn:5:Emu2}
	\end{align}
	If the integrand here is non-zero, we necessarily have $|v_3| \leq r$, $|v_2+v_3| \leq r$,
	and $|\lambda v_1+v_3|\geq R$. 
	This gives,
	\begin{equation}
		| \lambda v_1 | \geq | \lambda v_1 + v_3 | - |v_3| \geq \frac{R}{2} \;\;\text{ and }\;\;
		| v_2 | \leq | v_2 + v_3 | +  |v_3| \leq r.
		\label{eqn:5:impliedineqs}
	\end{equation}
	We will use these inequalities to impose constraints on $|\lambda v_2 + v_1|$ and $|\lambda v_2 - v_1|$,
		which are the inputs to the Gaussians in representation \eqref{eqn:5:EBL2}.
		If we can ensure that these are large, the fast decay of the Gaussians will imply that $\Ecub$ is small.
	By inspection, we see that large values of $|\lambda|$ pose a problem, but such large values can be dealt with separately
		by using the decay of $1/(1+\lambda^2)$ in \eqref{eqn:5:EBL2}.

	\emph{Regime One}: $|\lambda|\leq \sqrt{R}/4$.
	Observe that,
	$
		| \lambda v_2 + v_1 | \geq |v_1| - |\lambda v_2| 
                \geq R/(2|\lambda|) - 2|\lambda|,
	$
	in the last step using \eqref{eqn:5:impliedineqs}.
	Because the function $x \mapsto R/(2x) - 2x$ is decreasing for positive $x$, we have,
	\[
		| \lambda v_2 + v_1 | \geq \frac{R}{2|\lambda|} - 2|\lambda| \geq \frac{R}{2 (\sqrt{R}/4)} - 2(\sqrt{R}/4) > \sqrt{R}.
	\]
	An identical argument shows that $|\lambda v_2 - v_1| > \sqrt{R}$.
	It follows that,
	\begin{align*}
		A :&= \frac{1}{\sqrt{2} \pi} \int_{|\lambda| \leq \sqrt{R}/4} \frac{1}{1+\lambda^2} | E_{B(\lambda)}( G,f_1,f_2,G,f_3,f_4) | d\lambda	\\
			&= \frac{1}{\sqrt{2} \pi} \int_{|\lambda| \leq \sqrt{R}/4} \frac{1}{1+\lambda^2}| E_{B(\lambda)}( G \chi_{|x|\geq \sqrt{R}} ,f_1,f_2,G \chi_{|x| \geq \sqrt{R}},f_3,f_4) | d\lambda	\\
			&\leq \frac{1}{\sqrt{2} \pi } \left( \int_\R \frac{1}{1+\lambda^2} d\lambda \right)
				\| G \chi_{|x|\geq \sqrt{R} } \|_{L^2}^2 
				\| f_1 \|_{L^2}
				\| f_2 \|_{L^2}
				\| f_3 \|_{L^2}
				\| f_4 \|_{L^2}.
	\end{align*}
	We estimate,
	\[
		\| G \chi_{|x|\geq \sqrt{R} } \|_{L^2}^2 
			= 2\int_{\sqrt{R}}^\infty e^{-x^2} dx
			\leq \frac{2}{ \sqrt{R} } \int_{\sqrt{R}}^\infty x e^{-x^2} dx 
			= \frac{2}{ \sqrt{R} } e^{-R} \leq \frac{2}{ \sqrt{R} },
	\]
	and hence,
	$
		A \leq  (C/\sqrt{R})
				\| f_1 \|_{L^2}
				\| f_2 \|_{L^2}
				\| f_3 \|_{L^2}
				\| f_4 \|_{L^2},
	$
	for some absolute constant $C$.

	\emph{Regime Two}: $|\lambda|\geq \sqrt{R}/4$.
	This is easier: we have,
	\begin{align*}
		B :&= \frac{1}{\sqrt{2} \pi} \int_{|\lambda| \geq \sqrt{R}/4} \frac{1}{1+\lambda^2} | \E_{B(\lambda)}( G,f_1,f_2,G,f_3,f_4) | d\lambda	\\
			&\leq \frac{\sqrt{2}}{\pi}  \left( \int_{\sqrt{R}/4}^\infty \frac{1}{1+\lambda^2} d\lambda \right)
				\| G \|_{L^2}^2 
				\| f_1 \|_{L^2}
				\| f_2 \|_{L^2}
				\| f_3 \|_{L^2}
				\| f_4 \|_{L^2}.
	\end{align*}
	We estimate
	\[
		\int_{\sqrt{R}/4}^\infty \frac{1}{1+\lambda^2} d\lambda  \leq
		\int_{\sqrt{R}/4}^\infty \frac{1}{\lambda^2} d\lambda  \leq = \frac{4}{\sqrt{R}},
	\]
	which then gives,
	$
		B \leq  (C/\sqrt{R})
				\| f_1 \|_{L^2}
				\| f_2 \|_{L^2}
				\| f_3 \|_{L^2}
				\| f_4 \|_{L^2}.
	$

	Then, because
		$\Ecub(f_1,f_2,f_3,f_4) = A + B$, equation \eqref{eqn:5:refined1} is established.
\end{proof}

\begin{thm}
	There is an absolute constant $C$ such that if $\mu\in [0,1/2)$,
		$f_k$ is supported in $B(0,r)$ and $f_j$ is supported in $B(0,R)^C$, with $R>4r$, then
	\begin{equation}
		| \Ecub^\mu(f_1,f_2,f_3,f_4) |\leq 
			\frac{C}{(1-2\mu)^{5/8} R^{1/4} }
			\| f_1 \|_{L^2}
			\| f_2 \|_{L^2}
			\| f_3 \|_{L^2}
			\| f_4 \|_{L^2}
			\label{eqn:5:refined}
	\end{equation}
\end{thm}

\begin{proof}
	For fixed $\mu \in [0,1/2)$, let $\alpha = \sqrt{(1/2)-\mu}$.
	We will again adopt the notation $f^\lambda(x) = \lambda^{1/2} f(\lambda x)$.
	With this notation we have,
	\[
		e^{-(\frac{1}{2}-\mu)x^2} = e^{-(\alpha x)^2} = G(\alpha x) = \alpha^{-1/2} G^\alpha(x).
	\]
	Using the scaling property of $E_{B(\lambda)}$, we have,
	\begin{align*}
		\Ecub^\mu(f_1,f_2,f_3,f_4) 
			&=\frac{1}{\sqrt{2}\pi} \int_\R \frac{1}{1+\lambda^2} E_{B(\lambda)}
			\left(
				\alpha^{-1/2} G^\alpha, f_1, f_2,
				\alpha^{-1/2} G^\alpha, f_3, f_4
			\right)		\\
			&=\frac{1}{\alpha} \frac{1}{\sqrt{2}\pi} \int_\R \frac{1}{1+\lambda^2} E_{B(\lambda)}
			\left(
				 G, f_1^{1/\alpha}, f_2^{1/\alpha},
				 G, f_3^{1/\alpha}, f_4^{1/\alpha}
			\right) 	\\
			&= \frac{1}{\alpha} \Ecub \left(  f_1^{1/\alpha}, f_2^{1/\alpha}, f_3^{1/\alpha}, f_4^{1/\alpha} \right).
	\end{align*}

	Now assume that $f_{k_1}$ is supported in $B(0,R)^C$ and $f_{k_2}$ is supported in $B(0,r)$.
	We then have that $f_{k_1}^{1/\alpha}$ is supported in $B(0,\alpha R)^C$ and $f_{k_2}^{1/\alpha}$ is supported in $B(0,\alpha r)$.
	Then, using representation \eqref{eqn:5:Ehermite}, we have,
	\begin{align*}
		|\Ecub^\mu(f_1,f_2,f_3,f_4)| 
		&= 
			\left| \frac{1}{\alpha} \Ecub(  f_1^{1/\alpha}, f_2^{1/\alpha}, f_3^{1/\alpha}, f_4^{1/\alpha} )	\right| \\
			&=
			\left| \frac{1}{\alpha} \frac{2}{\pi} \int_0^{\pi/2} \int_\R 
			(e^{itH}f_1^{1/\alpha})
			(e^{itH}f_2^{1/\alpha})
			\overline{ ( e^{itH}f_3^{1/\alpha} )(e^{itH}f_4^{1/\alpha}}) dx dt	\right| \\
			&\leq \frac{1}{\alpha} \left( 
			\frac{2}{\pi} \int_0^{\pi/2} \int_\R  |(e^{itH}f_{k_1}^{1/\alpha}) 		(e^{itH}f_{k_2}^{1/\alpha})|^2 \right)^{1/2}
			\left( \frac{2}{\pi} \int_0^{\pi/2} \int_\R  |(e^{itH}f_{k_3}^{1/\alpha}) 		(e^{itH}f_{k_4}^{1/\alpha})|^2 \right)^{1/2} \\
			&= \frac{1}{\alpha} 
				\Ecub(  f_{k_1}^{1/\alpha}, f_{k_2}^{1/\alpha}, f_{k_1}^{1/\alpha}, f_{k_2}^{1/\alpha} )^{1/2}	
				\Ecub(  f_{k_3}^{1/\alpha}, f_{k_4}^{1/\alpha}, f_{k_3}^{1/\alpha}, f_{k_4}^{1/\alpha} )^{1/2}.
	\end{align*}
	Using \eqref{eqn:5:refined1}, we get,
	\[
		\Ecub(  f_{k_1}^{1/\alpha}, f_{k_2}^{1/\alpha}, f_{k_1}^{1/\alpha}, f_{k_2}^{1/\alpha} )^{1/2}
			\leq \frac{C}{ (\alpha R)^{1/2} } 
			\|  f_{k_1}^{1/\alpha} \|_{L^2}^2 
			\|  f_{k_2}^{1/\alpha} \|_{L^2}^2 
			= \frac{C}{ (\alpha R)^{1/2} } 
			\|  f_{k_1} \|_{L^2}^2 
			\|  f_{k_2} \|_{L^2}^2 ,
	\]
	while for the other $\Ecub$ term we can use the usual $L^2$ boundedness.
	This gives,
	\[
		|\Ecub^\mu(f_1,f_2,f_3,f_4)|  \leq \frac{ C }{ \alpha^{5/4} R^{1/4} } \|f_1\|_{L^2} \|f_2\|_{L^2} \|f_3\|_{L^2} \|f_4\|_{L^2}.
	\]
	Substituting back in $\alpha = \sqrt{(1/2)-\mu}$ gives the result.
\end{proof}

\subsubsection{Stationary waves are analytic}

\begin{thm}
	Suppose that $\phi \in L^2$ is a stationary wave solution of $iu_t = \Tcub(u,u,u)$; that is,
	$\phi$ satisfies,
	\begin{equation}
		\omega \phi (x)= \Tcub(\phi,\phi,\phi)(x),
		\label{eqn:5:strongSW}
	\end{equation}
	for some $\omega$.
	Then there exists $\alpha>0$ and $\beta>0$ such that
	$\phi e^{\alpha x^2} \in L^\infty$ and $\widehat{\phi} e^{\beta x^2} \in L^\infty$.
	As a result, $\phi$ can be extended to an entire function on the complex plane.
	\label{thm:5:smooth}
\end{thm}

Using the proof of Corollary \ref{thm:4:smooth}, this theorem is an immediate consequence of the following proposition.
\begin{prop}
	Suppose that $\phi\in L^2$ satisfies
	\begin{equation}
		\omega |\phi(x)| \leq \Tcub( |\phi|, |\phi|, |\phi| )(x),
		\label{eqn:5:weakSW}
	\end{equation}
	for some $\omega >0$.
	Then there exists $\alpha>0$ such that $x\mapsto \phi(x) e^{\alpha x^2} \in L^2$.
\end{prop}

\begin{proof}[Proof of proposition]
	For the proof, we will find $\mu$ so that we have the bound
	$\| \phi G_{\mu,\epsilon} \|_{L^2} \lesssim 1$
	independently of $\epsilon$. 
	Taking the limit $\epsilon \rightarrow 0$ will the yield the result.
	The structure of proof here is extremely similar to that of Theorem \ref{thm:4:strongL2}.
	For brevity, we will only describe the start of the proof here, which is the only part that is essentially different to the proof of Theorem \ref{thm:4:strongL2}.

	First, we fix throughout $\mu\leq 1/4$.
	Using formulas \eqref{eqn:5:Emubound1}, \eqref{eqn:5:Emubound2} and \eqref{eqn:5:refined}, there are constants $C$ independent of $\mu$, such that,
	\begin{align}
		\Ecub(f_1,f_2,f_3,f_4 \Gmu) 	&\leq \Ecub^\mu(f_1\Gmu,f_2\Gmu,f_3\Gmu,f_4)	\label{eqn:5:Emufact1} \\
		|\Ecub^\mu(f_1,f_2,f_3,f_4)| 	&\leq  C \| f_1 \|_{L^2}\| f_2 \|_{L^2}\| f_3 \|_{L^2}   \| f_4 \|_{L^2} 	\label{eqn:5:Emufact2} \\
		|\Ecub^\mu(f_1,f_2,f_3,f_4)| 	&\leq  C \frac{1}{ R^{1/4} } \| f_1 \|_{L^2}\| f_2 \|_{L^2}\| f_3 \|_{L^2}   \| f_4 \|_{L^2}, \label{eqn:5:Emufact3}  
	\end{align}
	where in the last inequality, $f_i$ is supported in $B(0,r)$ and $f_j$ is supported in $B(0,R)^C$ for $R>4r$.

	Now consider a function $\phi$ satisfying \eqref{eqn:5:weakSW}.
	We may assume $\phi$ is non-negative.
	For any $M>0$ define,
	\begin{align*}
		\phi_<(x) &= \phi(x) \chi_{|x|\leq M}(x),
		&\phi_{\sim}(x) &= \phi(x) \chi_{M<|x|\leq M^2}(x),
		&\phi_>(x) &= \phi(x) \chi_{|x|\leq M^2}(x).
	\end{align*}
	We have the decomposition $\phi = \phi_> + \phi_{\sim} + \phi_>$, and the supports are all disjoint,  which gives,
	\[
		\| \phi \Gmu \|_{L^2}^2 =
		\| \phi_< \Gmu \|_{L^2}^2 +
		\| \phi_{\sim} \Gmu \|_{L^2}^2 +
		\| \phi_> \Gmu \|_{L^2}^2.
	\]

	The first two terms are trivial to to bound uniformly in $M$.
	If $|x|\leq M^2$, we have,
	\[
		\Gmu(x) \leq e^{\mu x^2} \leq e^{\mu M^4}
	\]
	so setting  $\mu = M^{-4}$ gives
	$
		\| \phi_< \Gmu \|_{L^2} \leq 
		\| \phi_< e^1 \|_{L^2} \leq e^1 \| \phi \|_{L^2} \lesssim 1$, with the same bound for $\phi_\sim$.
	In order to prove the theorem, it remains then to bound 
		$\| \phi_> e^{G_{\mu,\epsilon,n}} \|_{L^2} $.

	Starting with the equation \eqref{eqn:5:weakSW} of the theorem, we multiply both sides by $\phi_>(x)  \Gmu(x)^2 $
	which gives,
	\[
		\omega  \phi_>(x) ^2 \Gmu(x)^2  \leq \Tcub(\phi,\ldots,\phi)(x) \phi_>(x)\Gmu(x)^2 
	\]
	Now integrating over $\R$ and using \eqref{eqn:5:Emufact1} gives,
	\begin{align*}
		\omega \| \phi_> \Gmu  \|_{L^2}^2& \leq \Ecub(\phi,\phi,\phi,\phi_> \Gmu^2)	\leq\Ecub^\mu( \phi \Gmu,\phi \Gmu,\phi \Gmu, \phi_> \Gmu).
	\end{align*}
	For convenience, let $\psi = \phi \Gmu $.
	The bound then reads,
	\[
		\omega \| \psi_> \|_{L^2}^2 \lesssim \Ecub^\mu(\psi,\psi,\psi,\psi_>).
	\]
	Now write each $\psi = \psi_< + \psi_{\sim} + \psi_>$ and expand the multinear functional.
	We will get many terms, which we bound in one of two ways.

	\begin{itemize}
		\item If there are three or more $\psi_>$ terms, bound by $\|\psi_>\|_{L^2}^k$ where $k$ is the number of $\psi_>$ terms appearing, using \eqref{eqn:5:Emufact2}.
		In this case the other terms are $\psi_<$ or $\psi_\sim$, which are uniformly bounded.

		\item
		If there are one or two $\psi_>$ terms, then there is either a $\psi_<$ term or a $\psi_\sim$ term.
			In the former case we can use the refined multilinear estimate \eqref{eqn:5:Emufact3}, with $R=M^2$,
			and bound by $M^{-1/2} \|\psi_>\|^k$ 
		(where $k=1$ or $k=2$).
			In the latter case we can bound by $\| \psi_\sim \|_{L^2} \|\psi_> \|_{L^2}^k \lesssim \| \phi_\sim \|_{L^2} \|\psi_> \|_{L^2}^k$ using \eqref{eqn:5:Emufact1}.
	\end{itemize}

	In total, we get,
	\begin{align*}
		\omega \| \psi_> \|_{L^2}^2 &\leq  \Ecub^\mu(\psi,\ldots,\psi,\psi_>)	\\
			&\leq C \left(  \| \psi_> \|_{L^2}^4 + \| \psi_> \|_{L^2}^3 + (M^{-1/2} + \| \phi_\sim \|_{L^2}) (
				\|\psi_>\|_{L^2}^2 +
				\|\psi_>\|_{L^2} ) \right),
	\end{align*}
	for a constant $C$ independent of $\mu$.
	This formula has the same structure as equation \eqref{eqn:4:poly} in the proof of Theorem \ref{thm:4:strongL2}.
	Replicating the same argument there,
		we find that if we choose $M$ sufficiently large there is a constant independent of $\epsilon$ such that $\| \psi_> \|_{L^2} \leq C$.
	Letting $\epsilon \rightarrow 0 $ then gives the result.
\end{proof}

\section*{Acknowledgments}

I am indebted to my doctoral advisor, Pierre Germain, for suggesting the problems considered in this work,
	and for many extensive mathematical discussions.
My thanks also to Zaher Hani for his conversations and lectures during the 2016 Hamiltonian PDE summer school in Maiori, Italy, 
	and to Miles Wheeler at the Courant Institute.

\bibliographystyle{amsplain}
\bibliography{references}

\end{document}